\documentclass[12pt]{article}
\usepackage{setspace}
\usepackage{indentfirst}
\usepackage{float}
\usepackage{graphicx}
\usepackage{bm}
\usepackage{amsmath}
\usepackage{amsthm}
\usepackage{amssymb}
\usepackage{subcaption}

\usepackage{placeins}

\usepackage{diagbox}
\usepackage{makecell}
\usepackage{authblk}

\usepackage[margin=1in,footskip=0.25in]{geometry}

\usepackage[pdftex,colorlinks,
            urlcolor = blue,
            linkcolor=blue,
            anchorcolor=blue,
            citecolor=blue
            ]{hyperref}
\usepackage[round]{natbib}
            
\newtheorem{proposition}{Proposition}

\title{Gibbs Phenomenon Suppression in PDE-Based Statistical Spatio-Temporal Models}
\author{Guanzhou Wei$^a$, Xiao Liu$^{a,}$\thanks{Corresponding author. E-mail address: \href{mailto:xl027@uark.edu}{xl027@uark.edu}.}, Russell Barton$^b$} 

\affil{\small$^a$Department of Industrial Engineering, University of Arkansas, USA}
\affil{$^b$Smeal College of Business, The Pennsylvania State University, USA}
\date{}

\begin{document}
\maketitle

\begin{abstract}
A class of physics-informed spatio-temporal models has recently been proposed for modeling spatio-temporal processes governed by advection-diffusion equations. The central idea is to approximate the process by a truncated Fourier series and let the governing physics determine the dynamics of the spectral coefficients. 
However, because many spatio-temporal processes in real applications are non-periodic with boundary discontinuities, the well-known Gibbs phenomenon and ripple artifact almost always exist in the outputs generated by such models due to truncation of the Fourier series. 
Hence, the key contribution of this paper is to propose a physics-informed spatio-temporal modeling approach that significantly suppresses the Gibbs phenomenon when modeling spatio-temporal advection-diffusion processes.  The proposed approach starts with a data flipping procedure for the process respectively along the horizontal and vertical directions (as if we were unfolding a piece of paper that has been folded twice along the two directions). Because the flipped process becomes spatially periodic and has a complete waveform without any boundary discontinuities, the Gibbs phenomenon disappears even if the Fourier series is truncated. Then, for the flipped process and given the Partial Differential Equation (PDE) that governs the process, this paper extends an existing PDE-based spatio-temporal model by obtaining the new temporal dynamics of the spectral coefficients, while maintaining the physical interpretation of the flipped process. Numerical investigations based on a real dataset have been performed to demonstrate the advantages of the proposed approach. It is found that the proposed approach effectively suppresses the Gibbs Phenomenon and significantly reduces the ripple artifact in modeling spatio-temporal advection-diffusion processes. 
Computer code is available on GitHub. 
\end{abstract}
\noindent\textbf{Key words:} {\em spatial-temporal models, advection-diffusion processes, physics-informed statistical learning, Gibbs phenomenon suppression, non-periodic processes, boundary discontinuities.}
\clearpage

\setstretch{1.25}
\section{Introduction} \label{sec:intro}

\subsection{Background and Motivating Examples}
Statistical modeling of spatio-temporal data has been extensively investigated in physical science, engineering, and environment, e.g., the modeling of temperature fields \citep{guinness2013interpolation,kuusela2018locally,vandeskog2022quantile}, polutant propagation \citep{huang2004modeling,liu2016statistical,schliep2020data,fioravanti2022spatiotemporal}, wildfires, precipitation and extreme weather events \citep{cooley2007bayesian,heaton2011spatio,kleiber2012daily,liu2018spatio,bopp2021hierarchical,wei2022}, biomedical and remote-sensing images \citep{kang2011meta,katzfuss2011spatio,hefley2017dynamic,reich2018fully,castruccio2018scalable,zammit2021deep}, sensor placement and spatial design \citep{zimmerman2006optimal,zimmerman2019environmental,sharrock2022joint} etc. In particular, for spatio-temporal data arising from advection-diffusion processes, a class of physics-informed statistical modeling approach has been proposed which represents the spatial process by its (truncated) Fourier series and let the governing physics (usually Partial Differential Equations (PDE)) determine the stochastic temporal dynamics of the spectral coefficients \citep{Cressie2011, Sigrist2015, Kutz2016, Liu2021}. Such an approach is based on the classical solution of nonlinear dynamical systems using the spectral theory and eigenfunction expansions \citep{Kutz2016}.

\begin{figure}[h!]  
	\begin{center}
		\includegraphics[width=\textwidth]{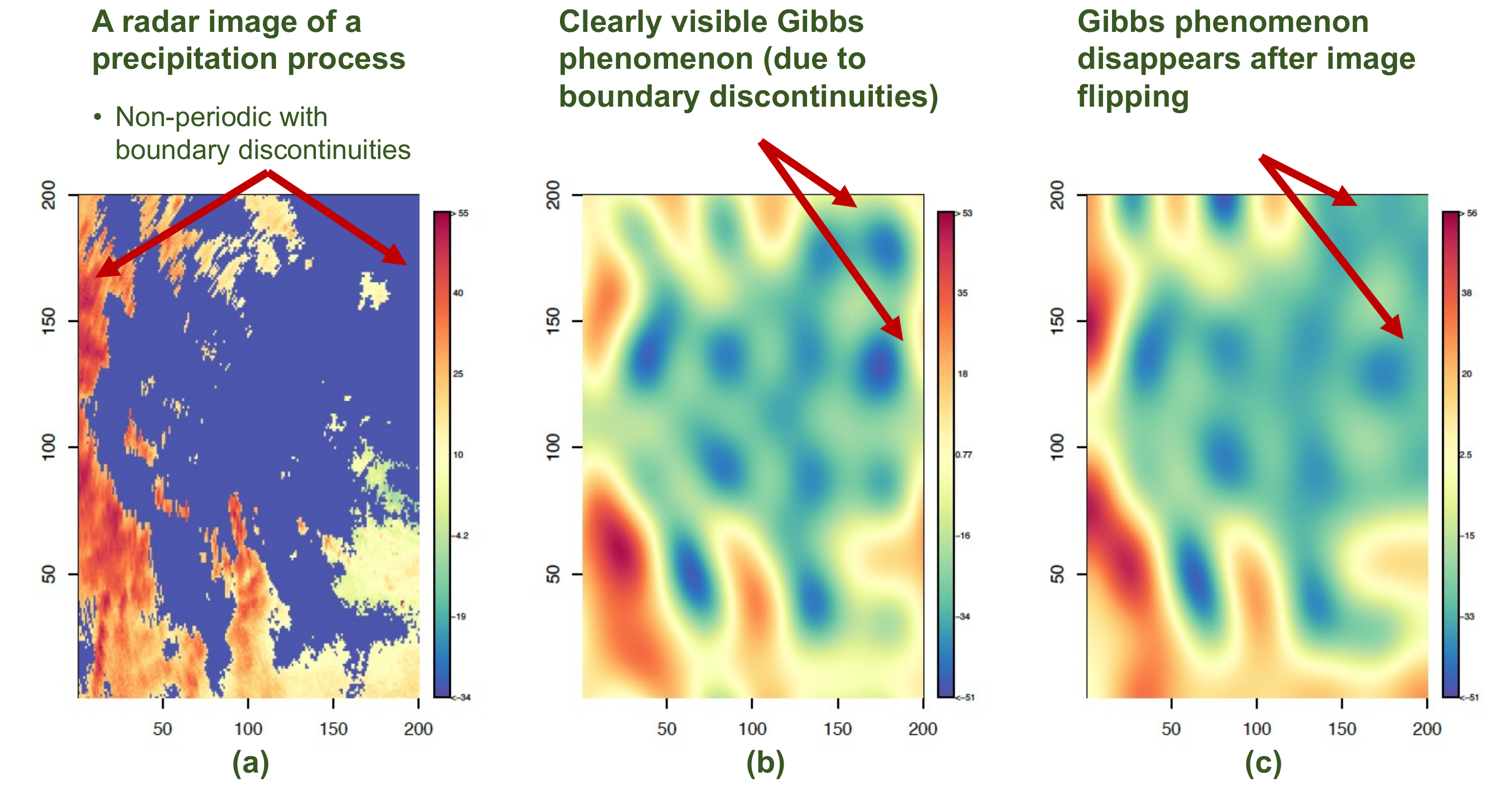}
		\centering
		\vspace{-6pt}
		\caption{Illustration of the Gibbs Phenomenon: (a) weather radar reflectivity field of a precipitation process with discontinuous boundaries; (b) a low-pass filtered image re-constructed from a truncated Fourier series; (c) a low-pass filtered image re-constructed by a truncated Fourier series after image flipping (details are given in Section \ref{sec:flip}).}
		\label{fig:illustration}
	\end{center}
\end{figure}

However, the Gibbs Phenomenon is a well-known issue when truncating a Fourier series \citep{Pan1993, Jerri1998}. In spatial data modeling, such a phenomenon is caused by the discontinuities of the spatial process at the boundaries of the spatial domain. Although the Fourier series requires the spatial process to be periodic, many real spatio-temporal processes  do not satisfy this condition. For illustrative purposes, Figure \ref{fig:illustration}(a) shows a weather radar reflectivity field of a precipitation process. It is clearly seen that this process is non-periodic with boundary discontinuities; for example, the values on the left and right boundaries of the spatial domain do not match. As a result, the Gibbs phenomenon is clearly visible in Figure \ref{fig:illustration}(b), which shows a reconstructed image from a truncated Fourier series (i.e., the higher-order Fourier modes are discarded). Near the top-right boundaries of Figure \ref{fig:illustration}(b), we could clearly observe positive radar intensity, while in effect there should be no positive readings in that region based on Figure \ref{fig:illustration}(a). This is known as the Gibbs Phenomenon: when there exist boundary discontinuities, the truncated Fourier series does not have a sufficient number of terms to handle the abrupt jump of the radar intensity value between the two boundaries (discontinuities). We leave the explanation of Figure \ref{fig:illustration}(c) to Section \ref{sec:flip}. Of course, when the spatial process does not have strong boundary discontinuities, the Gibbs Phenomenon may not occur. For example, when the weather system is relatively small and is only located in the central area of the radar image, the boundary values are all or close to zero and no strong Gibbs Phenomenon is expected.

\subsection{Literature Review and Contributions}
Existing methods to suppress the Gibbs phenomenon include windowing, data flipping, baseline tilting, digital comparison, etc. Most of these methods aim to complete the periodicity of the waveform in one way or another. For example, the well-known Hanning or Hamming method uses window functions to force the endpoints of the waveform to be zero or near zero \citep{Prabhu2014}. The baseline tilting approach achieves periodicity by adjusting each data point to force the last point of the waveform to be equal to the first point \citep{Pan2001}. The digital comparison method applies a filter that chooses the truncation length such that the endpoint of that truncated signal is the same or close to the first point of the truncated signal.  The former two approaches inevitably distort the original data and lead to less accurate models (imagine applying a window function on the spatial domain to force the boundaries values to be equal), while the latter approach cannot be applied to spatial data where the spatial domain is fixed and cannot be truncated. 

Among the existing methods, one approach appears to be particularly promising for our problems---the data flipping approach. This approach was initially proposed by \cite{Pan1993} for one-dimensional time-dependent signals. The idea is to flip the signal (i.e., to create the mirror image of the data) such that the flipped signal is completely periodic (with equal start and end points). This technique creates a new signal, which is twice the length of the original signal, that completely removes the discontinuity at the endpoints and removes the Gibbs phenomenon. 

In our problem, on the other hand, we are dealing with spatial data at a given time. Hence, the data flipping technique, originally developed for one-dimensional time-varying signals, needs to be extended. Then, for the new flipped process, the physics-informed or PDE-based spatio-temporal model proposed in \cite{Liu2021} also needs to be extended to handle the flipped spatio-temporal data and retain the physical interpretations of the model. Hence, this paper makes the following contributions. 

$\bullet$ The paper firstly proposes the data flipping approach for spatial data. Unlike the data flipping for one-dimensional time-varying signals, it is necessary to flip the spatial data twice, respectively along the horizontal and vertical directions (as if we were unfolding a piece of paper that has been folded twice along the two directions). The flipped image has the complete waveform and the Gibbs phenomenon disappears; see Figure \ref{fig:illustration}(c). Details are provided in Section \ref{sec:flip}.  

$\bullet$ Based on the double flipped data, this paper extends the PDE-based spatio-temporal model described in \cite{Liu2021}. In that paper, the authors obtained the  temporal dynamics of the spectral coefficients according to the governing advection-diffusion PDE. In this paper, we obtain the temporal dynamics of the spectral coefficients for the flipped process, while maintaining the physical interpretation of such temporal dynamics according to the governing advection-diffusion PDE. Details are provided in Section \ref{sec:model}. 

$\bullet$ Numerical investigations, based on both simulated and real data sets, have been performed to demonstrate the advantages of the proposed approach over existing methods. Findings are provided in Section \ref{sec:example}. 

\section{Gibbs Phenomenon Suppression in PDE-Based Statistical Spatio-Temporal Models} \label{sec:one}

\subsection{Data Flipping} \label{sec:flip}
\begin{figure}[ht!]  
	\begin{center}
		\includegraphics[width=\textwidth]{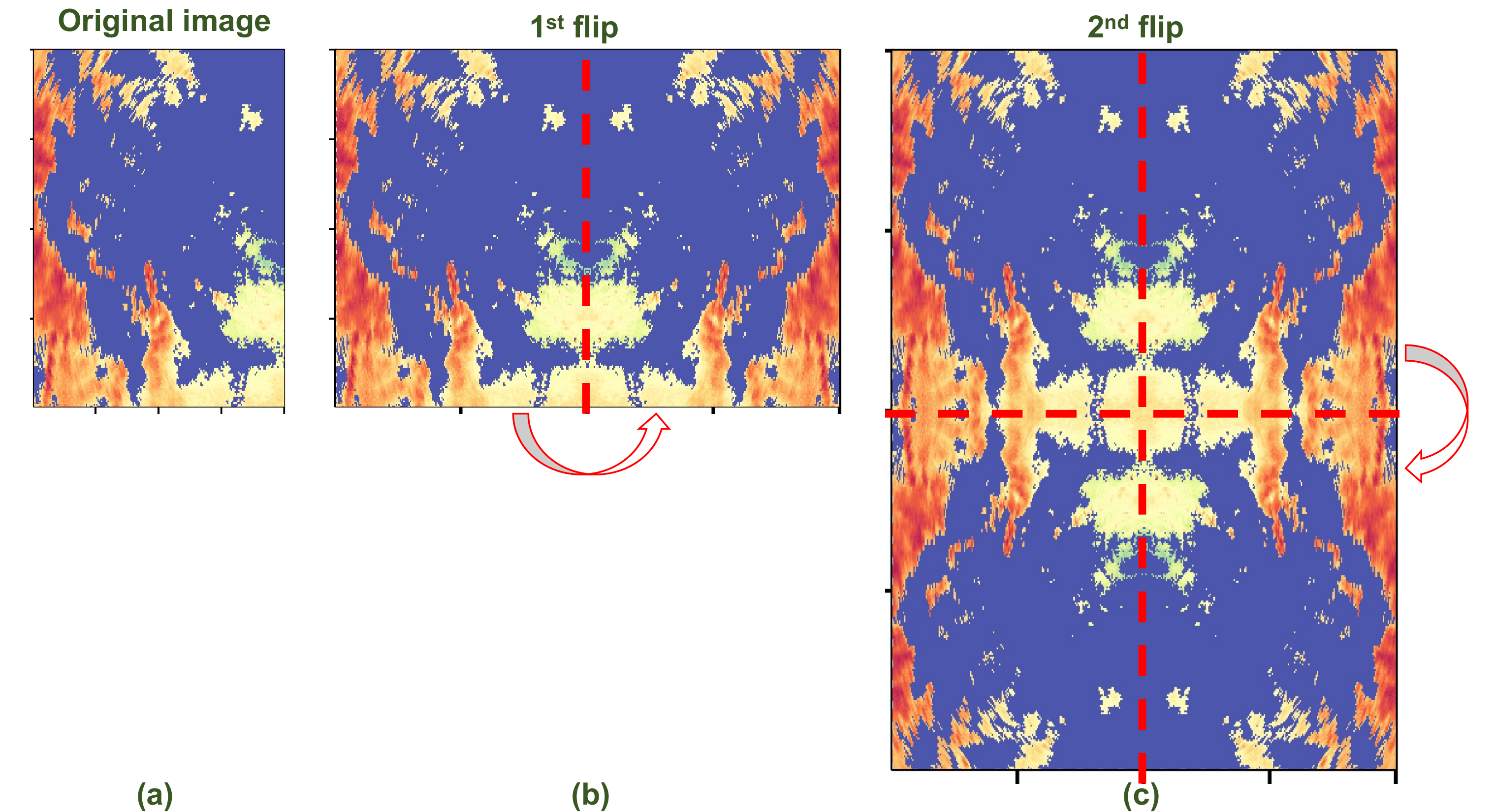}
		\centering
		\vspace{-12pt}
		\caption{Flip the spatial data twice along the $x$- and $y$-direction to create a periodic process}
		\label{fig:flip}
	\end{center}
\end{figure}

Because the Gibbs phenomenon arises from the incomplete periodicity of the spatial image, flipping is an effective technique that makes the boundary values of a spatial process to be equal and thus creates the complete periodicity. The flipping process is illustrated in Figure \ref{fig:flip}. Unlike the data flipping for one-dimensional time-dependent signals, two flips are needed for spatial data. The first flip creates the mirror image of the original image along the $x$-direction (Figure \ref{fig:flip}(b)), and the second flip creates the mirror image of the flipped image along the $y$-direction (Figure \ref{fig:flip}(c)). Note that, it is also possible to flip the image along the $y$-direction first and then the $x$-direction. 

As seen in Figure \ref{fig:flip}(c), the flipped image has complete periodicity (with equal boundary values). If we re-visit Figure \ref{fig:illustration}(c) in Section 1, it is seen that the Gibbs phenomenon disappears completely in the image reconstructed by the truncated Fourier series. Unlike the existing windowing methods to eliminate the Gibbs phenomenon, the flipping technique does not force weights to distort the original data. Also note that, the flipped image is four times larger than the original image, and the spectral resolution is also four times higher, potentially increasing the computational burden of the algorithm. This is the necessary cost in order to suppress the Gibbs phenomenon.



To formally define the flipping process, consider a spatial-temporal process defined on a $N_1\times N_2$ uniform spatial grid given by a Cartesian product $\mathbb{S} = \mathcal{X} \times \mathcal{Y} = \{\bm{s}=(x, y)|x\in \mathcal{X}, y\in \mathcal{Y}\}$, where $\mathcal{X} = \{0,1/N_1,\dots, (N_1-1)/N_1\}$ and $\mathcal{Y}=\{0, 1/N_2,\dots,(N_2-1)/N_2\}$. Let $N= N_1\times N_2$, and $\bm{y} = (y_1, y_2, \cdots, y_{N})^T$ be a vector of all observations from the spatial grid with $y_i=y(\bm{s}_i)$ being the observation at location $\bm{s}_i$, $i=1,2,\cdots,N$.

As shown in Figure \ref{fig:flip}, the flipped image is defined on a $2N_1\times 2N_2$ uniform spatial grid. If we let $\bm{y}^* = (y^*_1, y^*_2, \cdots, y^*_{4N})^T$ be a vector that contains all observations from the flipped image, it is possible to show that
\begin{equation} \label{eq:rotation}
\bm{y}^*=\bm{R}\bm{y}
\end{equation}
where $\bm{R}$ is a $4N \times N$ \textit{flipping matrix} given by
\begin{equation} \label{eq:flipping}
    \bm{R} = 
               [ \bm{I}_{N_2} , \bm{J}_{N_2} ]^T \otimes [\bm{I}_{N_1}, \bm{J}_{N_1}]^T
\end{equation}
with $\bm{I}$ and $\bm{J}$ respectively being an identity matrix and the mirror image of an identity matrix (i.e., with 1s on the bottom left to upper right diagonal line). 

\begin{figure}[ht!]  
	\begin{center}
		\includegraphics[width=0.7\textwidth]{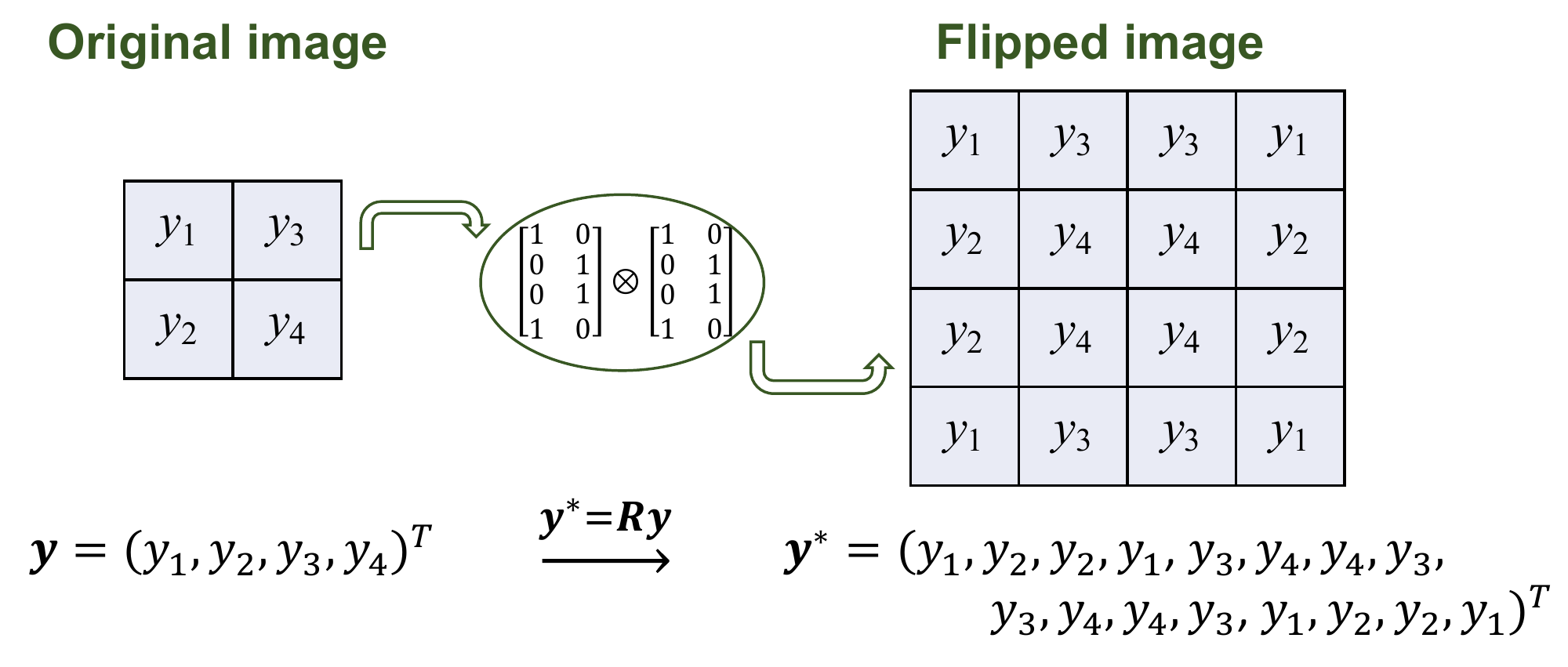}
		\centering
		\vspace{-10pt}
		\caption{An illustration of how the observation vector of the flipped image is created through the flipping operation} 
		\label{fig:flip2}
	\end{center}
\end{figure}

For an illustrative example, consider a simple $2\times2$ image with $\bm{y} = (y_1, y_2, y_3, y_4)^T$ in Figure \ref{fig:flip2}. Given the rotation matrix
\begin{equation*}
\bm{R} = \begin{bmatrix}
            1 & 0 & 0 & 1  \\
 	 0 & 1 & 1& 0
         \end{bmatrix}^T \otimes 
         \begin{bmatrix}
           1 & 0 & 0 & 1  \\
 	 0 & 1 & 1& 0
         \end{bmatrix}^T, 
\end{equation*}
the operation $\bm{y}^* = \bm{R} \bm{y}$ creates a new vector which is exactly the vector by stacking the columns of the flipped image. 
\vspace{-6pt}


To show why (\ref{eq:rotation}) holds, consider a $N_2\times N_1$ pixel array $\bm{M}$ of the original image. Then, the corresponding pixel array $\bm{M}^*$ of the flipped image can be obtained as:
\begin{equation}
    \bm{M}^* = [\bm{I}_{N_2} ,  \bm{J}_{N_2}]^T \bm{M} [\bm{I}_{N_1}, \bm{J}_{N_1}].
\end{equation}

Here, it is noted that, multiplying a permutation matrix $[\bm{I}_{N_1}, \bm{J}_{N_1}]$ to the right of $\bm{M}$ yields the first flip along the $x$-direction, and multiplying $[\bm{I}_{N_2}, \bm{J}_{N_2}]^T$ to the left yields the second flip along the $y$-direction. Hence, let $\text{vec}(\bm{M}^*)\triangleq\bm{y}^*$ and $\text{vec}(\bm{M})\triangleq\bm{y}$, we have
\begin{equation}
    \begin{split}
        \bm{y}^* & \triangleq \text{vec}(\bm{M}^*)\\
                 & = \text{vec}([\bm{I}_{N_2}, \bm{J}_{N_2}]^T\bm{M}[\bm{I}_{N_1}, \bm{J}_{N_1}])\\
                 & =  [ \bm{I}_{N_2} , \bm{J}_{N_2} ]^T \otimes [\bm{I}_{N_1}, \bm{J}_{N_1}]^T\text{vec}(M) \\
                 & \triangleq \bm{Ry}
    \end{split}
\end{equation}

Note that, the places of $\bm{I}_{N_2}$ and $\bm{J}_{N_2}$ in (\ref{eq:flipping}) are interchangeable, and so are the places of $\bm{I}_{N_1}$ and $\bm{J}_{N_1}$. The four place settings of $\bm{I}_{N_1}$ and $\bm{J}_{N_1}$, and $\bm{I}_{N_2}$,  $\bm{J}_{N_2}$ generate four different flipping matrices that control the flipping sequence along the $x$-direction and the $y$-direction. For example, although the setting in (\ref{eq:rotation}) firstly flips the image along the $x$-direction using the right boundary, and then flips the flipped image along the $y$-direction using the bottom boundary, any of the four flipping matrices is a valid choice that generates a periodic process with no boundary discontinuities.

\subsection{The Model} \label{sec:model}

For the flipped process, this section presents a PDE-based spatio-temporal model to suppress the Gibbs phenomenon, while maintaining the physical interpretation of the model according to the governing advection-diffusion PDE. 

\subsubsection{Preliminaries}
Consider a real-valued advection-diffusion process given by a PDE:
\begin{equation}\label{eq:advection-diffusion} 
\dot{\xi}(\bm{s},t) = \mathcal{A}\xi(\bm{s},t)+ Q(\bm{s},t), 
\end{equation} 
where $\xi(\bm{s},t)$ is the spatio-temporal process for some physical quantities, $Q(\bm{s},t)$ is the forcing term, and $\mathcal{A}$ is an advection-diffusion operator: 
\begin{equation}\label{eq:operator} 
	\mathcal{A}\xi(\bm{s},t) = -\bm{v}^T(\bm{s},t)\nabla \xi(\bm{s},t)+\nabla\cdot[\bm{D}(\bm{s},t)\nabla \xi(\bm{s},t)],
\end{equation} 
where $\bm{v}(\bm{s}, t)$, $\bm{D}(\bm{s}, t)$, $\nabla$, and $\nabla\cdot$ are the velocity, diffusivity, gradient, and divergence, respectively. 


By projecting the process $\xi(\bm{s},t)$ and the external forcing $Q(\bm{s}, t)$ onto a subspace spanned by Fourier bases, the processes $\xi(\bm{s}, t)$ and $Q(\bm{s}, t)$ can be expressed by a linear combination of orthogonal spatial basis functions and temporal coefficients that evolve over time \citep{wikle1998hierarchical, Sigrist2015, Liu2021}:
\begin{equation} \label{eq:decomposition2}
	\xi(\bm{s},t) = \sum_{\bm{k}\in\mathcal{K}_1}\alpha_{\bm{k}}^{(c)}(t)f_{\bm{k}}^{(c)}(\bm{s}) + 2\sum_{\bm{k}\in\mathcal{K}_2} \left(\alpha_{\bm{k}}^{(c)}(t)f_{\bm{k}}^{(c)}(\bm{s}) + \alpha_{\bm{k}}^{(s)}(t)f_{\bm{k}}^{(s)}(\bm{s})\right),
\end{equation}
\begin{equation} \label{eq:decompositionQ}
	Q(\bm{s}, t) =\sum_{\bm{k}\in\mathcal{K}_1}\beta_{\bm{k}}^{(c)}(t)f_{\bm{k}}^{(c)}(\bm{s})+2\sum_{\bm{k}\in\mathcal{K}_2} ({\beta}_{\bm{k}}^{(c)}(t)f_{\bm{k}}^{(c)}(\bm{s}) + \beta_{\bm{k}}^{(s)}(t)f_{\bm{k}}^{(s)}(\bm{s}))
\end{equation}
where $f_{\bm{k}}^{(c)}(\bm{s})=\cos(\bm{k}^T\bm{s})$ and $f_{\bm{k}}^{(s)}(\bm{s})=\sin(\bm{k}^T\bm{s})$, $\alpha_{\bm{k}}^{(c)}(t)$, $\alpha_{\bm{k}}^{(s)}(t)$, $\beta_{\bm{k}}^{(c)}(t)$, and $\beta_{\bm{k}}^{(s)}(t)$ are the temporal coefficients, and $\bm{k}$ is the spatial wavenumber taken from the sets:
\begin{equation*}\label{eq:K1K2}
	\begin{split}
	\mathcal{K}_1 &= \left\{(0,0), \left(0, N_2\big/2\right), \left(N_1\big/2, 0\right), \left(N_1\big/2, N_2\big/2\right)\right\}\\
\mathcal{K}_2 &= \left\{(k_1,k_2)|k_1=1,\cdots, N_1\big/2-1, k_2=-N_2\big/2+1,\cdots, -1\right\}\\
	&\hskip0.2in\cup\left\{(k_1,k_2)|k_1=0,1,\cdots, N_1\big/2, k_2=1,2,\cdots, N_2\big/2-1\right\}\\
	&\hskip0.2in\cup\left\{\left(k_1,N_2\big/2\right)|k_1=0,1,\cdots, N_1\big/2\right\}.
	\end{split}
\end{equation*}

Because the Fourier basis is pre-determined, the modeling of the spatio-temporal process can be converted to the modeling of the temporal stochastic processes of $\alpha_{\bm{k}}^{(c)}(t)$ and $\alpha_{\bm{k}}^{(s)}(t)$. Following this idea, a PDE-based spatio-temporal model can be written as follows \cite{Liu2021}:
\begin{equation} \label{eq:model}
\begin{split}
    \bm{y}(t) & = \bm{F}\bm{\alpha}(t) + \text{noise}\\
    \bm{\dot{\alpha}}(t) & = \bm{P}(t)\bm{\alpha}(t) + \bm{\beta}(t) + \text{noise}
    \end{split}
\end{equation}
where 

$\bullet$ $\bm{y}(t) = (y(\bm{s}_1, t), y(\bm{s}_2, t), \cdots, y(\bm{s}_N, t))^T$ is a vector that contains all observations at time $t$ from spatial locations $\bm{s}_1, \bm{s}_2, \cdots, \bm{s}_N$. 

$\bullet$  $\bm{F}=(\bm{f}(\bm{s}_1),\cdots,\bm{f}(\bm{s}_{N}))^T\in\mathbb{R}^{N\times (2|\mathcal{K}_1|+|\mathcal{K}_2|)}$ contains the Fourier basis functions, where $\bm{f}(\bm{s})=(\bm{f}_{\mathcal{K}_1}^{(c)}(\bm{s}), 2\bm{f}_{\mathcal{K}_1}^{(c)}(\bm{s}), 2\bm{f}_{\mathcal{K}_2}^{(c)}(\bm{s}))^T$, $\bm{f}_{\mathcal{K}_1}^{(c)}(\bm{s})=(f^{(c)}_{\bm{k}_{1,1}}(\bm{s}),f^{(c)}_{\bm{k}_{1,2}}(\bm{s})\cdots,f^{(c)}_{\bm{k}_{1,|\mathcal{K}_1|}}(\bm{s}))$, $\bm{f}_{\mathcal{K}_2}^{(c)}(\bm{s})=(f^{(c)}_{\bm{k}_{2,1}}(\bm{s}),f^{(c)}_{\bm{k}_{2,2}}(\bm{s}),\cdots,f^{(c)}_{\bm{k}_{2,|\mathcal{K}_2|}}(\bm{s}))$, and $\bm{f}_{\mathcal{K}_2}^{(s)}(\bm{s})=(f^{(s)}_{\bm{k}_{2,1}}(\bm{s}),f^{(s)}_{\bm{k}_{2,2}}(\bm{s}),\cdots,f^{(s)}_{\bm{k}_{2,|\mathcal{K}_2|}}(\bm{s}))$.

$\bullet$ $\bm{\alpha}(t)$ and $\bm{\beta}(t)$ are the spectral coefficients from (\ref{eq:decomposition2}) and (\ref{eq:decompositionQ}). In particular, let $\bm{k}_{1,i}$ be the $i$-th component in set $\mathcal{K}_1$, and $\bm{k}_{2,i}$ be the $i$-th component in set $\mathcal{K}_2$, then, $\bm{\alpha}(t)=\text{vec}(\bm{\alpha}^{(c)}_{\mathcal{K}_1}(t), \bm{\alpha}^{(c)}_{\mathcal{K}_2}(t), \bm{\alpha}^{(s)}_{\mathcal{K}_2}(t))$, where 
$\bm{\alpha}^{(c)}_{\mathcal{K}_1}(t)=(\alpha^{(c)}_{\bm{k}_{1,1}}(t), \alpha^{(c)}_{\bm{k}_{1,2}}(t),\cdots, \alpha^{(c)}_{\bm{k}_{1,{|\mathcal{K}_1|}}}(t))^T$,
$\bm{\alpha}^{(c)}_{\mathcal{K}_2}(t)=(\alpha^{(c)}_{\bm{k}_{2,1}}(t), \alpha^{(c)}_{\bm{k}_{2,2}}(t),\cdots, \alpha^{(c)}_{\bm{k}_{2,{|\mathcal{K}_2|}}}(t))^T$, and 
$\bm{\alpha}^{(s)}_{\mathcal{K}_2}(t)=(\alpha^{(s)}_{\bm{k}_{2,1}}(t), \alpha^{(s)}_{\bm{k}_{2,2}}(t),\cdots, \alpha^{(s)}_{\bm{k}_{2,{|\mathcal{K}_2|}}}(t))^T$. Similarly, $\bm{\beta}(t)=\text{vec}(\tilde{\bm{\beta}}^{(c)}_{\mathcal{K}_1}(t), \tilde{\bm{\beta}}^{(c)}_{\mathcal{K}_2}(t), \tilde{\bm{\beta}}^{(s)}_{\mathcal{K}_2}(t))$, where 
$\tilde{\bm{\beta}}^{(c)}_{\mathcal{K}_1}(t)=(\tilde{{\beta}}^{(c)}_{\bm{k}_{1,1}}(t), \tilde{{\beta}}^{(c)}_{\bm{k}_{1,2}}(t), \cdots, \tilde{{\beta}}^{(c)}_{\bm{k}_{1,{|\mathcal{K}_1|}}}(t))^T$,
$\tilde{\bm{\beta}}^{(c)}_{\mathcal{K}_2}(t)=(\tilde{{\beta}}^{(c)}_{\bm{k}_{2,1}}(t), \tilde{{\beta}}^{(c)}_{\bm{k}_{2,2}}(t), \cdots, \tilde{\beta}^{(c)}_{\bm{k}_{2,{|\mathcal{K}_2|}}}(t))^T$, and 
$\tilde{\bm{\beta}}^{(s)}_{\mathcal{K}_2}(t)=(\tilde{\beta}^{(s)}_{\bm{k}_{2,1}}(t),\tilde{\beta}^{(s)}_{\bm{k}_{2,2}}(t),\cdots, \tilde{\beta}^{(s)}_{\bm{k}_{2,{|\mathcal{K}_2|}}}(t))^T$.

$\bullet$ The transition matrix $\bm{P}(t)$ is found by the Galerkin method and is given by 
\begin{equation} \label{eq:P}
\bm{P}(t) = \begin{pmatrix}
	\bm{C}_1^{-1}\bm{\Psi}^{(\text{A1+D1})}_{11}& 2\bm{C}_{1}^{-1}\bm{\Psi}^{(\text{A1+D1})}_{12}& 2\bm{C}_{1}^{-1}\bm{\Psi}^{(\text{A2+D2})}_{12}\\
	\frac{1}{2}\bm{C}_{2}^{-1}\bm{\Psi}^{(\text{A1+D1})}_{21}& 
	\bm{C}_{2}^{-1}\bm{\Psi}^{(\text{A1+D1})}_{22}&
	\bm{C}_{2}^{-1}\bm{\Psi}^{(\text{A2+D2})}_{22}\\
	\frac{1}{2}\bm{C}_{2}^{-1}\bm{\Psi}^{(\text{A3+D3})}_{21}& 
	\bm{C}_{2}^{-1}\bm{\Psi}^{(\text{A3+D3})}_{22}& 
	\bm{C}_{2}^{-1}\bm{\Psi}^{(\text{A4+D4})}_{22}
\end{pmatrix}.
\end{equation}

Detailed expressions of the transition matrix (\ref{eq:P}) can be found in Appendix A. It is critical to note that, the transition matrix is directly derived from the underlying advection-diffusion process, and depends on important physical parameters such as the advection velocity and diffusivity. Hence, the PDE-based spatio-temporal model (\ref{eq:model}) has a strong physical interpretation and integrates governing physics into a data-driven model. We refer the readers to \cite{Sigrist2015, Liu2021, wei2022} for more details. 


\subsubsection{The Extended Model for the Flipped Process}

In model (\ref{eq:model}), the dimension of $\bm{\alpha}(t)$ is $N = N_1\times N_2$. For large $N$, one commonly used approach is to truncate $\bm{\alpha}(t)$ by retaining only the low-frequency components. As discussed in Section \ref{sec:intro}, when the spatio-temporal process is not periodic with boundary discontinuities, removing the higher-order Fourier modes gives rise to the Gibbs phenomenon (Figure \ref{fig:illustration}(b)). To suppress the Gibbs phenomenon, we flip the process at any time $t$ using the approach described in Section \ref{sec:flip} and directly model the flipped process. Because the flipped process $\xi^*(\bm{s}, t)$ is periodic and has no discontinuities on the boundaries, truncating the spectral coefficients $\bm{\alpha}^*(t)$ of the flipped process no longer causes the Gibbs phenomenon. 

A key question immediately arises: \textit{what is the temporal dynamics of the spectral coefficients $\bm{\alpha}^*(t)$ of the flipped process $\xi^*(\bm{s}, t)$ according to the governing advection-diffusion equation}? Since the flipped process is obtained by flipping the original process governed by an advection-diffusion PDE, the temporal dynamics of $\bm{\alpha}^*(t)$ must also reflect the underlying governing physics. The answer to this question is given by the following proposition:

\begin{proposition}
Let $\xi(\bm{s}, t)$ be an advection-diffusion process governed by the PDE (\ref{eq:advection-diffusion}), and let $\xi^*(\bm{s}, t)$ be a flipped process obtained from (\ref{eq:rotation}). The temporal dynamics of the spectral coefficients $\bm{{\alpha}^*}(t)$ of the flipped process is given by an Ordinary Differential Equation (ODE) as follows:
\begin{equation}\label{eq:proposition1}
    \bm{\dot{\alpha}}^*(t) = \bm{H}\bm{P}(t)\bm{H}^\dag\bm{\alpha}^*(t) + \bm{\beta}^*(t)
\end{equation}
where $\bm{P}(t)$ is the transition matrix of $\bm{{\alpha}}(t)$, defined in (\ref{eq:P}), for the original process $\xi(\bm{s}, t)$, $\bm{H}=(\bm{F}^*)^{\dag}\bm{R}\bm{F}$, $\bm{F}^*$ contains the Fourier bases for the flipped process, $\bm{R}$ is the flipping matrix defined in (\ref{eq:rotation}), $\bm{\beta}^*(t)=\bm{H}\bm{\beta}(t)$, and $^\dag$ is the Moore-Penrose pseudo-inverse.
\end{proposition}

To see why Proposition 1 holds, note that $(\xi^*(\bm{s}_1,t),\cdots,\xi^*(\bm{s}_{4N},t))^T = \bm{R} (\xi(\bm{s}_1,t),\cdots,\xi(\bm{s}_N,t))^T$. Applying the Fourier transform on both sides, we obtain:
\begin{equation} \label{eq:temporalcoeffODE2}
    \bm{\alpha}^*(t)=(\bm{F}^*)^{\dag}\bm{R}\bm{F}\bm{\alpha}(t),
\end{equation}
where $\bm{F}$ and $\bm{F}^*$ are respectively the Fourier bases for $\xi(\bm{s}, t)$ and $\xi^*(\bm{s}, t)$.

Taking the derivative on both sides of ($\ref{eq:temporalcoeffODE2}$) yields:
\begin{equation} \label{eq:temporalcoeffODE3}
    \begin{split}
        \bm{\dot{\alpha}}^*(t) &= (\bm{F}^*)^{\dag}\bm{R}\bm{F}\bm{\dot{\alpha}}(t)\\
                               &= (\bm{F}^*)^{\dag}\bm{R}\bm{F}[\bm{P}(t)\bm{\alpha}(t) + \bm{\beta}(t)]\\
                               &= (\bm{F}^*)^{\dag}\bm{R}\bm{F}\bm{P}(t)\bm{F}^{\dag}\bm{R}^{-1}\bm{F}^*\bm{\alpha}^*(t) + (\bm{F}^*)^{\dag}\bm{R}\bm{F}\bm{\beta}(t)\\
                               &\triangleq \bm{H}\bm{P}(t)\bm{H}^{\dag}\bm{\dot{\alpha}}^*(t) + \bm{\beta}^*(t).
    \end{split}
\end{equation}
 
Comparing (\ref{eq:model}) and (\ref{eq:temporalcoeffODE3}), we see that $\bm{\beta}^*(t)=\bm{H}\bm{\beta}(t)=(\bm{F}^*)^{\dag}\bm{R}\bm{F}\bm{\beta}(t)$, which is equivalent to applying the Fourier transform on both the forcing term $\bm{Q}(\bm{s}, t)$ and the forcing term $\bm{Q}^*(\bm{s}, t)$ of the flipped process. 

Proposition 1 shows that the transition matrix of temporal coefficients $\bm{{\alpha}^*}(t)$ (of the flipped process) can be directly obtained from the transition matrix of $\bm{{\alpha}}(t)$ (of the original process) through $\bm{H}\bm{P}(t)\bm{H}^\dag$. Because the expression of $\bm{P}(t)$ has already been obtained in \cite{Liu2021} from the advection-diffusion PDE, one may easily compute $\bm{H}\bm{P}(t)\bm{H}^\dag$, which is also determined by the governing physics and thus possess a strong physical interpretation. 
If it is costly to compute the $N_1\times N_2$ matrix $\bm{P}(t)$, one often needs to truncate the coefficient $\bm{\alpha}^*(t)$ of the flipped process so that only a small number of components in $\bm{P}(t)$ are computed. Note that, truncating $\bm{\alpha}^*(t)$ of the flipped process (with no boundary discontinuities) no longer gives rise to the Gibbs phenomenon. 


Once the dynamics of temporal coefficients $\bm{\alpha}^*(t)$ has been established in (\ref{eq:temporalcoeffODE3}), we obtain a physics-informed statistical spatio-temporal model for the flipped process:
\begin{equation}\label{eq:proposition2}
	\begin{split}
	\bm{y}^*(t) &= (\bm{F}^*, \bm{0})\bm{\theta}^*(t) + \bm{v}^*(t), \hskip0.42in \bm{v}^*(t)\sim \mathcal{N}(\textbf{0},\bm{V}^*) \\
	\bm{\theta}^*(t)  &= \bm{G}^*(t)\bm{\theta}^*(t-\Delta) + \bm{w}^*(t), \hskip0.15in \bm{w}^*(t)\sim \mathcal{N}(\textbf{0},\bm{W}^*),
	\end{split}
\end{equation}
\vspace{-10pt}
where  $\bm{\beta}^*(t)$ is a Brownian motion with an unknown constant mean, and 
\begin{equation}
\bm{\theta}^*(t)=
\begin{pmatrix}
	\bm{\alpha}^*(t)\\
	\bm{\beta}^*(t)
\end{pmatrix},
\bm{G}^*(t)=
\begin{pmatrix}
	\exp(\Delta\bm{H}\bm{P}(t)\bm{H}^\dag) & \bm{I}\\
	\textbf{0} & \bm{I}
\end{pmatrix},
\bm{w}^*(t)=
\begin{pmatrix}
	\bm{w}_{\bm{\alpha}^*}(t)\\
	\bm{w}_{\bm{\beta}^*}(t)
\end{pmatrix}.
\end{equation}

The state transition equation in (\ref{eq:proposition2}) is obtained by solving (\ref{eq:temporalcoeffODE3}) using the first-order Euler forward method and the state augmentation approach \citep{evensen2000ensemble, stroud2007sequential, stroud2010ensemble}. Here, $\bm{w}_{\bm{\alpha}^*}(t)$ captures the uncertainty associated with the state transition as well as the numerical approximation error due to the first-order Euler forward method, and $\bm{w}_{\bm{\beta}^*}(t)$ captures the uncertainty associated with the Brownian motion $\bm{\beta}^*(t)$. It is important to mention that when the low-rank $\bm{P}(t)$ is used for computing $\bm{G}^*(t)$, biases are inevitably introduced into the transition of the spectral coefficients $\bm{\alpha}^*(t)$ due to the truncation of high-order Fourier modes. Hence, a bias correction process  can be added to the state transition equation in (\ref{eq:proposition2}) to capture the truncation error. However, the bias correction and the temporal coefficients $\bm{\beta}^*(t)$ of the forcing term  are not identifiable without prior knowledge or additional regularization. In such a case, the estimated $\bm{\beta}^*(t)$ captures both the forcing term and the bias correction process.

\section{Numerical Examples} \label{sec:example}
Two numerical examples are provided to illustrate the effectiveness of the proposed approach in suppressing the Gibbs phenomenon. The first example is based on a simulated dataset from  an advection process (Section \ref{sec:SimulationStudy}), and the second example is based on a real weather radar dataset for a spatio-temporal environmental process (Section  \ref{sec:CaseStudy}).
\vspace{-10pt}
\subsection{Example I}\label{sec:SimulationStudy}
In Example I, we demonstrate the capability of suppressing the Gibbs phenomenon using the proposed approach based on a simulated dataset from an advection process. In particular, we model the simulated datasets using both the approach described in \cite{Liu2021} (without data flipping) and the approach proposed in this paper (with data flipping). The performance of the two approaches is compared. 

The data used in Example I are simulated from the following advection process:
\begin{equation}\label{eq:CasePDE}
\begin{split}
    \dot{\xi}(\bm{s},t) &= -\bm{v}^T\nabla \xi(\bm{s},t) + Q(\bm{s})\\
    \xi(\bm{s},0) &=  Q(\bm{s})\\
    Q(\bm{s}) &= \frac{3}{2\pi\cdot0.18^2}\exp\left(-\frac{\|(0.1,0)-\bm{s}\|^2}{2\cdot0.18^2}\right)\\
     \bm{v} &= (0.01, 0)^T.
\end{split}
\end{equation} 

\begin{figure}[h!]
    \centering
    \begin{subfigure}[t]{0.5\linewidth}
        \centering
        \includegraphics[scale=0.65]{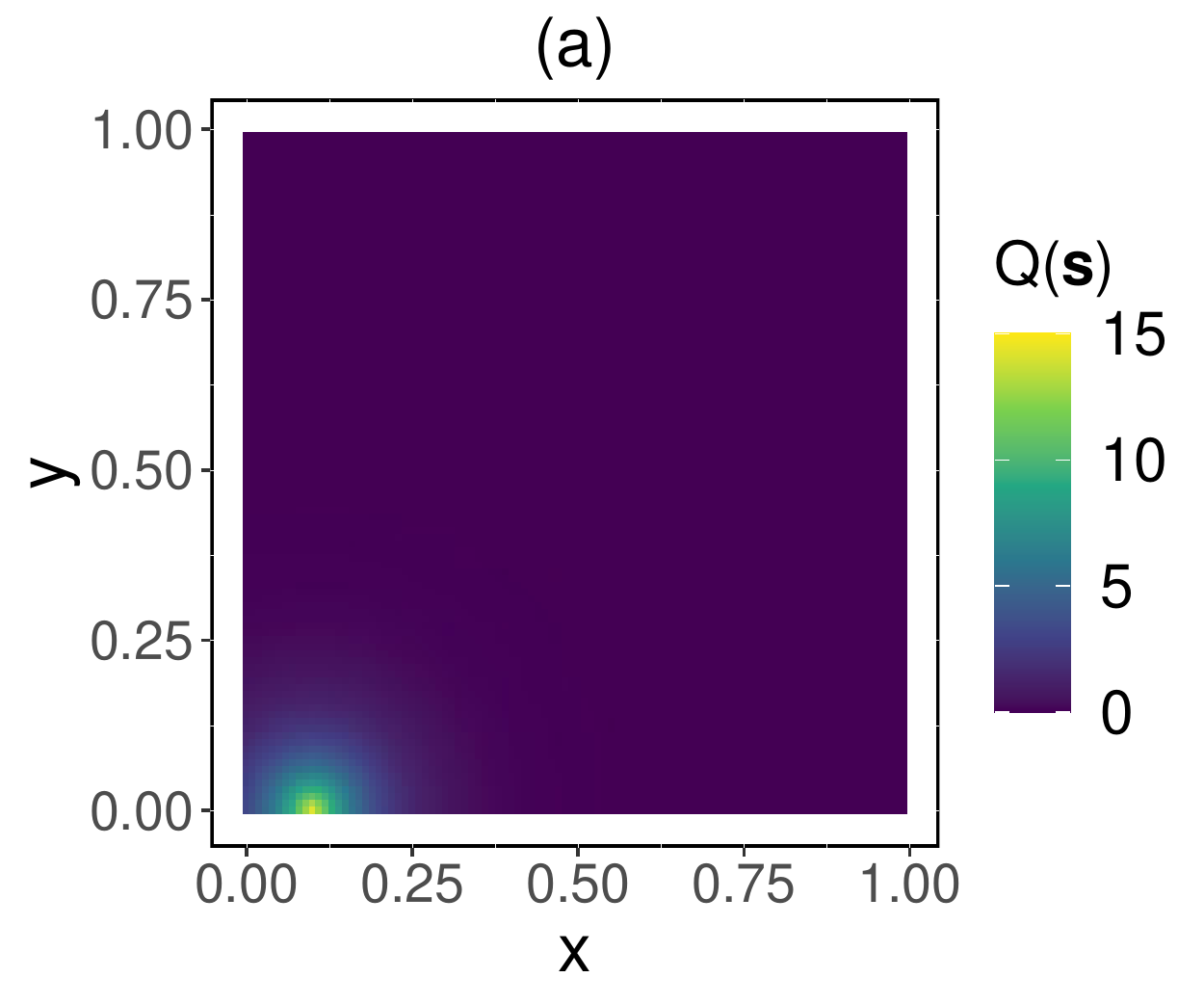}
    \end{subfigure}
    \quad
    \begin{subfigure}[t]{0.4\linewidth}
        \centering
        \includegraphics[scale=0.65]{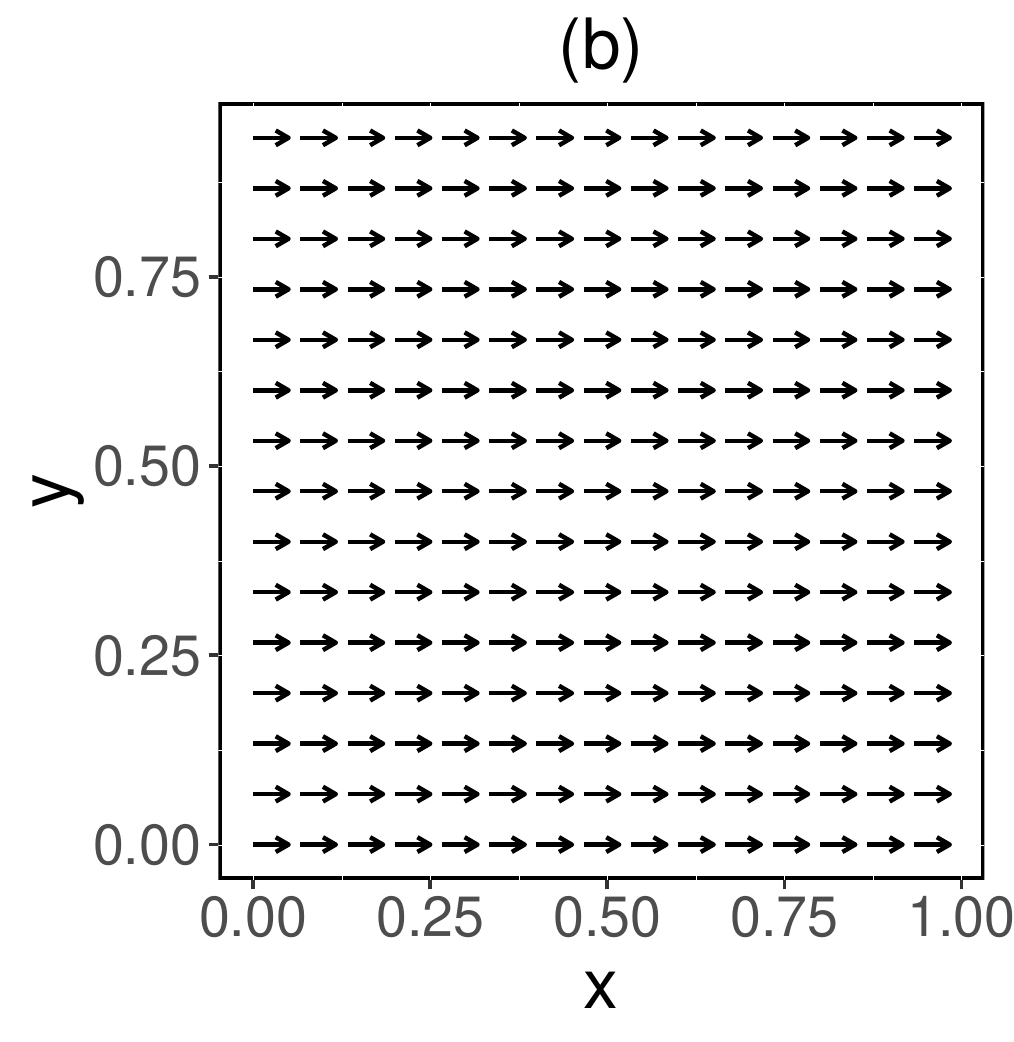}
	\end{subfigure}
	\caption{(a) The forcing term $Q(\bm{s})$; (b) the velocity field}\label{fig:SimulationSetting}
\end{figure}

The PDE (\ref{eq:CasePDE}) defines an advection process with a forcing term. In this example, we let the process be defined on a 100$\times$100 spatial grid over $[0,0.99]^2$. The forcing term is given by a function $Q(\bm{s})$ that is temporally invariant but spatially variant. Figure \ref{fig:SimulationSetting}(a) shows the image of the forcing term (which is also the initial condition for the process). We see that a signal source is placed near the bottom left corner of the spatial domain, which gives rise to a pronounced boundary discontinuity. The velocity field is chosen to be $\bm{v}= (0.01, 0)^T$ per time step, which does not change in time and space; see Figure \ref{fig:SimulationSetting}(b).

Solving the PDE of (\ref{eq:CasePDE}) using the spectral method, we obtain the dynamics of its temporal coefficients $\bm{\alpha}(t)$ in discrete time steps as:  
\begin{equation}\label{eq:CaseTemporalDynamics}
    \begin{split}
        \bm{\alpha}(t+\Delta) &= \exp(\Delta \bm{P}(t))\bm{\alpha}(t) + \bm{\beta}(t)\\
                \bm{\beta}(t) &= \bm{F}^{-1}Q(\bm{s})\\
                \bm{\alpha}(0) &= \bm{F}^{-1}Q(\bm{s}),
    \end{split}
\end{equation}
where $\bm{F}$ is defined below (\ref{eq:model}) and $\bm{P}(t)$ is given by (\ref{eq:P}). 

\begin{figure}[H]
    \centering
    \includegraphics[width=1\textwidth]{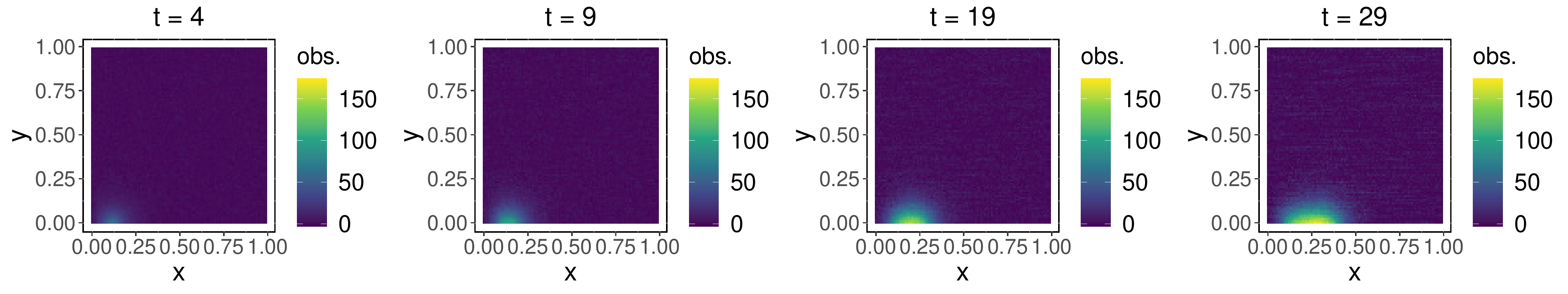}
    \caption{Simulated advection process for Example I}\label{fig:SimulatedProcess}
\end{figure}

To simulate noisy observations, two white noises with variances 0.005 and 0.001 are respectively added to $\bm{\alpha}(t)$ and $\bm{\beta}(t)$ in (\ref{eq:CaseTemporalDynamics}). By recursively computing the temporal coefficients $\bm{\alpha}(t)$ with (\ref{eq:CaseTemporalDynamics}), the simulated advection process with noisy observations is obtained after performing the Inverse Fourier Transform at discrete time steps. In particular, we simulate the process for 30 time steps at $t=0,1,2,\cdots,29$, where the first 20 images are used for constructing the model, and the last 10 images are used for investigating the predictive capabilities of the models. Figure \ref{fig:SimulatedProcess} shows the simulated advection process. 

Next, we model the simulated process using both the approach described in \cite{Liu2021} (without data flipping) and the approach proposed in this paper (with data flipping), and compare the performance of the two approaches. 

$\bullet$ (Existing Approach) We first model the simulated data using the existing physics-informed spatial-temporal model without data flipping. Following (\ref{eq:model}), the simulated data are modeled by a dynamical model:
\begin{equation}\label{eq:CaseDynamicalModel}
	\begin{split}
	    \bm{y}(t) &= (\bm{F},\bm{0})\bm{\theta}(t) + \bm{v}(t), \hskip0.42in \bm{v}(t)\sim\mathcal{N}(\bm{0}, \bm{V})\\
	    \bm{\theta}(t) &= \bm{G}\bm{\theta}(t-\Delta) + \bm{w}(t), \hskip0.32in \bm{w}(t)\sim\mathcal{N}(\bm{0}, \bm{W}), 
	\end{split}
\end{equation}
where $\bm{\theta}(t)= (\bm{\alpha}(t)^T,\bm{\beta}(t)^T)^T$,
$\bm{w}(t)=(\bm{w}_{\bm{\alpha}}^T(t),\bm{w}_{\bm{\beta}}^T(t))^T$, and
$\bm{G}(t)=
\begin{pmatrix}
	\exp(\Delta\bm{P}(t)) & \bm{I}\\
	\textbf{0} & \bm{I}
\end{pmatrix}$.

\begin{figure}[H]
    \centering
    \includegraphics[width=0.95\textwidth]{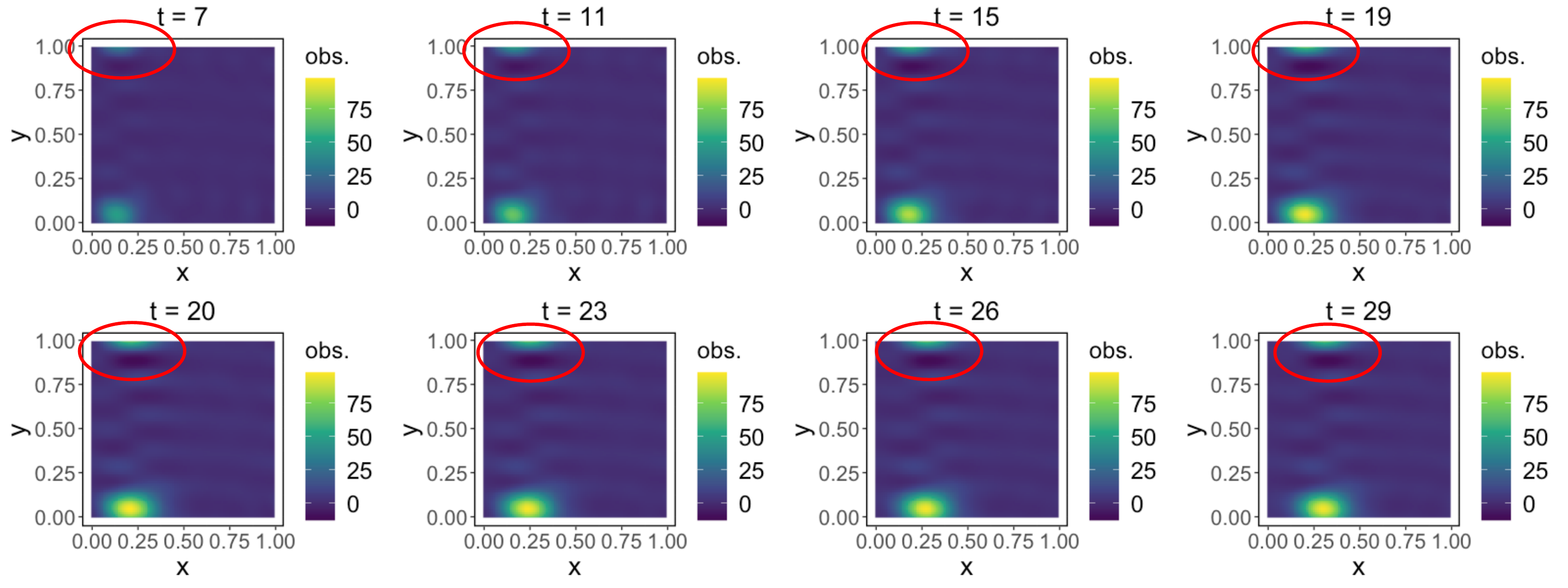}
    \caption{The filtered (top row) and predicted (bottom row) processes at different times, where both the modeling and prediction performance of the existing method may suffer from the Gibbs Phenomenon.}\label{fig:EstimatedSP}
\end{figure}

Even in this relatively simple example, the dimension of $\bm{\theta}(t)$ goes up to $2N_1\times N_2 = 20,000$, where $N_1=N_2=100$. Hence, a common practice is to truncate $\bm{\theta}(t)$ by retaining only the low-frequency portion of $\bm{\theta}(t)$. In this example, we set $ K = 100$, meaning that only 100 components in $\bm{\theta}(t)$ corresponding to the lower-order Fourier modes are kept. 
Then, the estimated (truncated) $\hat{\bm{\theta}}(t)$ is obtained by the Kalman Filter using the data from the first 20 time steps, and the approximated process is recovered from the Inverse Fourier Transform. The top row of Figure \ref{fig:EstimatedSP} shows the approximated process at times 7, 11, 15 and 19. The Gibbs Phenomenon is \textit{clearly visible} near the top left corner of the images (indicated by red circles). We clearly notice an area with elevated values, while the actual values in that area, shown in Figure \ref{fig:SimulatedProcess}, are almost zero. 
The bottom row of Figure \ref{fig:EstimatedSP} shows the predicted process at times 20, 23, 26 and 29. We again clearly observe the Gibbs Phenomenon (i.e., the area with elevated values indicated by red circles) from the predicted images. This observation shows that: \textit{both the modeling and prediction performance of the existing method may suffer from the Gibbs Phenomenon}. The cause for such a phenomenon is due to the boundary discontinuity (in this case, there exists a source term at the bottom boundary). In many real applications, such boundary discontinuities are widely observed. 

\begin{figure}[!ht]
    \centering
    \includegraphics[width=1\textwidth]{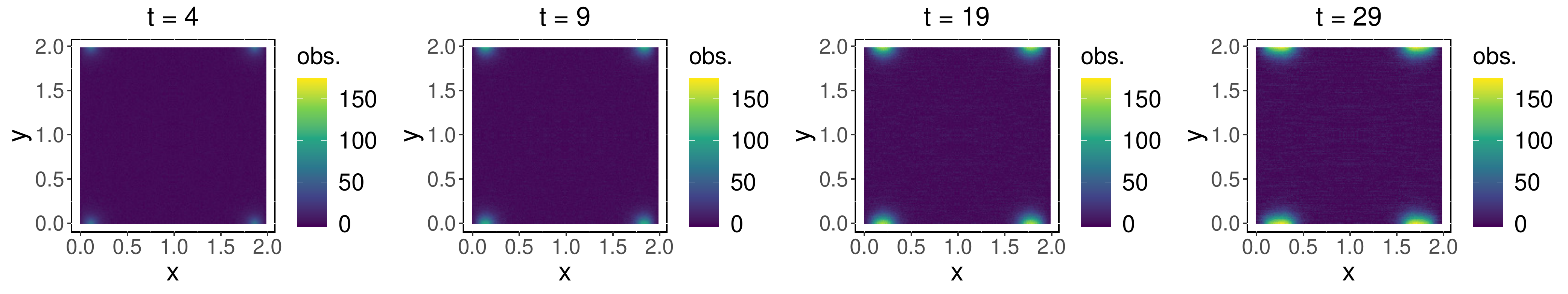}
    \caption{Flipped images of the advection process at times 4, 9, 19 and 29 with no boundary discontinuities.}\label{fig:FlippedSP}
\end{figure}

$\bullet$ (The Proposed Approach) Next, we model the same dataset using the proposed approach to suppress the Gibbs Phenomenon. As described in Section \ref{sec:flip}, the image at any time is firstly flipped twice to create a new image that is four times bigger than the original one. Figure \ref{fig:FlippedSP} shows the flipped images at times 4, 9, 19 and 29. It is seen that the new process becomes periodic and new images no longer have boundary discontinuities.

\quad

Based on (\ref{eq:proposition2}), we obtain the statistical model for the flipped images as follows:
\begin{equation}\label{eq:FSPDynamical}
    \begin{split}
        \bm{\tilde{y}}^*(t) &= (\bm{\tilde{F}^*}, \bm{0})\bm{\theta}^*(t) + \bm{\tilde{v}^*}(t) \hskip0.67in \bm{\tilde{v}}^*(t)\sim\mathcal{N}(\bm{0}, \bm{\tilde{V}}^*)\\
        \bm{\theta}^*(t) &= \bm{G}^*(t)\bm{\theta}^*(t-\Delta) + \bm{w}^*(t), \hskip0.32in \bm{w}^*(t)\sim\mathcal{N}(\bm{0}, \bm{W}^*),
    \end{split}
\end{equation}
where $\bm{\tilde{F}}^*=\bm{I}_{K^*}$, $\bm{\tilde{V}}^*=\bm{H}\sigma^2_{\bm{\alpha}}\bm{I}_{K^*}\bm{H}^T$, $\bm{W}^* = \begin{pmatrix}
	\bm{H}\sigma^2_{\bm{\alpha}}\bm{I}_{K^*}\bm{H}^T & \bm{0}\\
	\bm{0} & \bm{H}\sigma^2_{\bm{\beta}}\bm{I}_{K^*}\bm{H}^T
\end{pmatrix}$, and $K^*$ is the number of retained low-frequency coefficients for the flipped process. 

\begin{figure}[h!]
    \centering
    \includegraphics[width=1\textwidth]{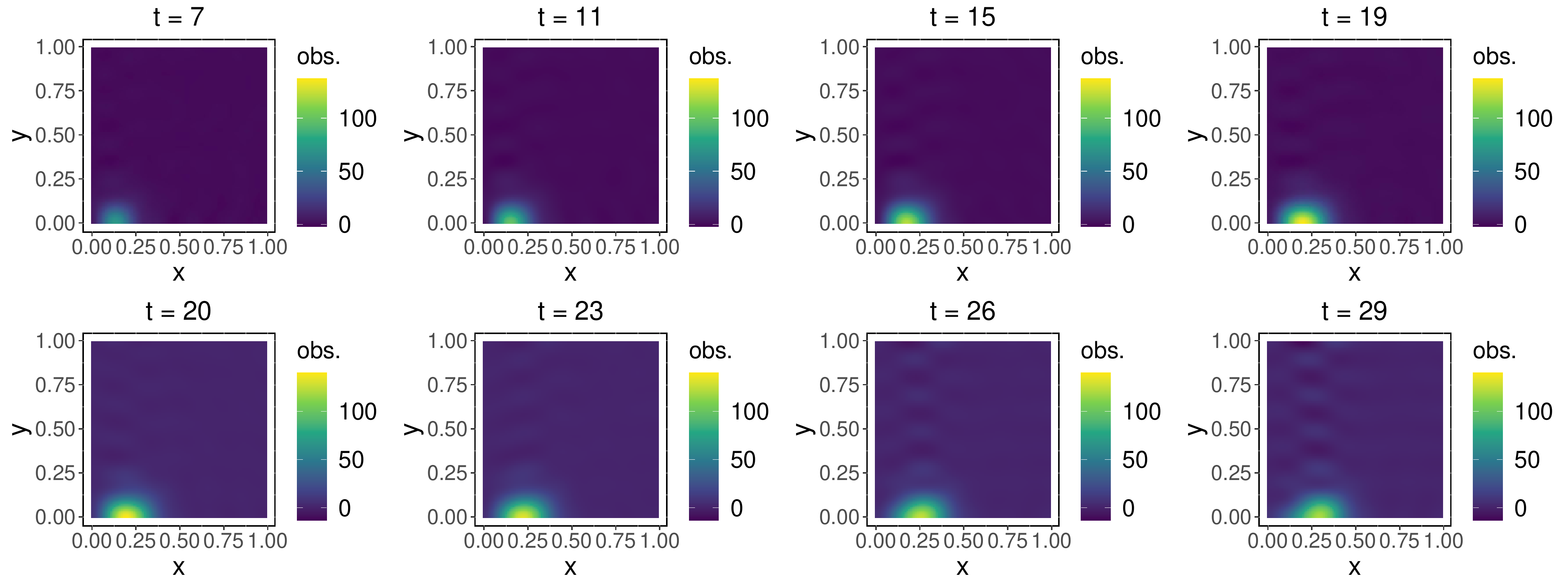}
    \caption{The filtered (top row) and predicted (bottom row) process using the proposed approach with the Gibbs Phenomenon being completely suppressed.}\label{fig:EstimatedFSP}
\end{figure}

Note that, the flipped images have a dimension of $2N_1\times 2N_2$, which is four times larger than the dimension of the original images. 
Hence, the number of low-frequency coefficients retained, $K^*=400$, is also chosen to be four times larger than $K=100$ when the existing method is used in the previous example. 
The Kalman Filter is performed to estimate the spectral coefficients using the data from the first 20 time steps, and the approximated process is recovered from the Inverse Fourier Transform. The top row of Figure \ref{fig:EstimatedFSP} shows the filtered process at times 7, 11, 15 and 19, while the bottom row of Figure (\ref{fig:EstimatedFSP}) shows the predicted process at  times 20, 23, 26 and 29. By comparing Figure \ref{fig:EstimatedSP} and Figure \ref{fig:EstimatedFSP}, it is immediately seen that \textit{the Gibbs Phenomenon is completed removed using the proposed approach}. In fact, it is also seen in Figure \ref{fig:EstimatedSP} that both the filtered and predicted images have some ``ripple 
artifact'' due to the truncation of the spectral coefficients. \textit{Such artifact has also been significantly reduced} by the proposed approach that constructs the model based on the flipped process. 


$\bullet$ (Comparison of Modeling and Prediction Accuracy) Removing the Gibbs Phenomenon and the ripple artifact potentially leads to a significant improvement of modeling accuracy. Using the same dataset, we compare the modeling accuracy between the proposed approach and the existing approach with different choices of $K$ (i.e., the number of low-frequency spectral coefficients retained). In particular, we construct three models using the existing approach by letting $K=100, 196, 400$ respectively, and compare the Mean-Absolute-Error (MAE) of the three models to the MAE of the model constructed using the proposed approach with $K^*=400$. We denote these four models by ``NF100", ``NF196", ``NF400", ``F400" with ``NF'' and ``F''  respectively denoting the existing approach with no process flipping and the proposed approach with process flipping.

\begin{figure}[h!]
    \centering
    \includegraphics[width=0.65\textwidth]{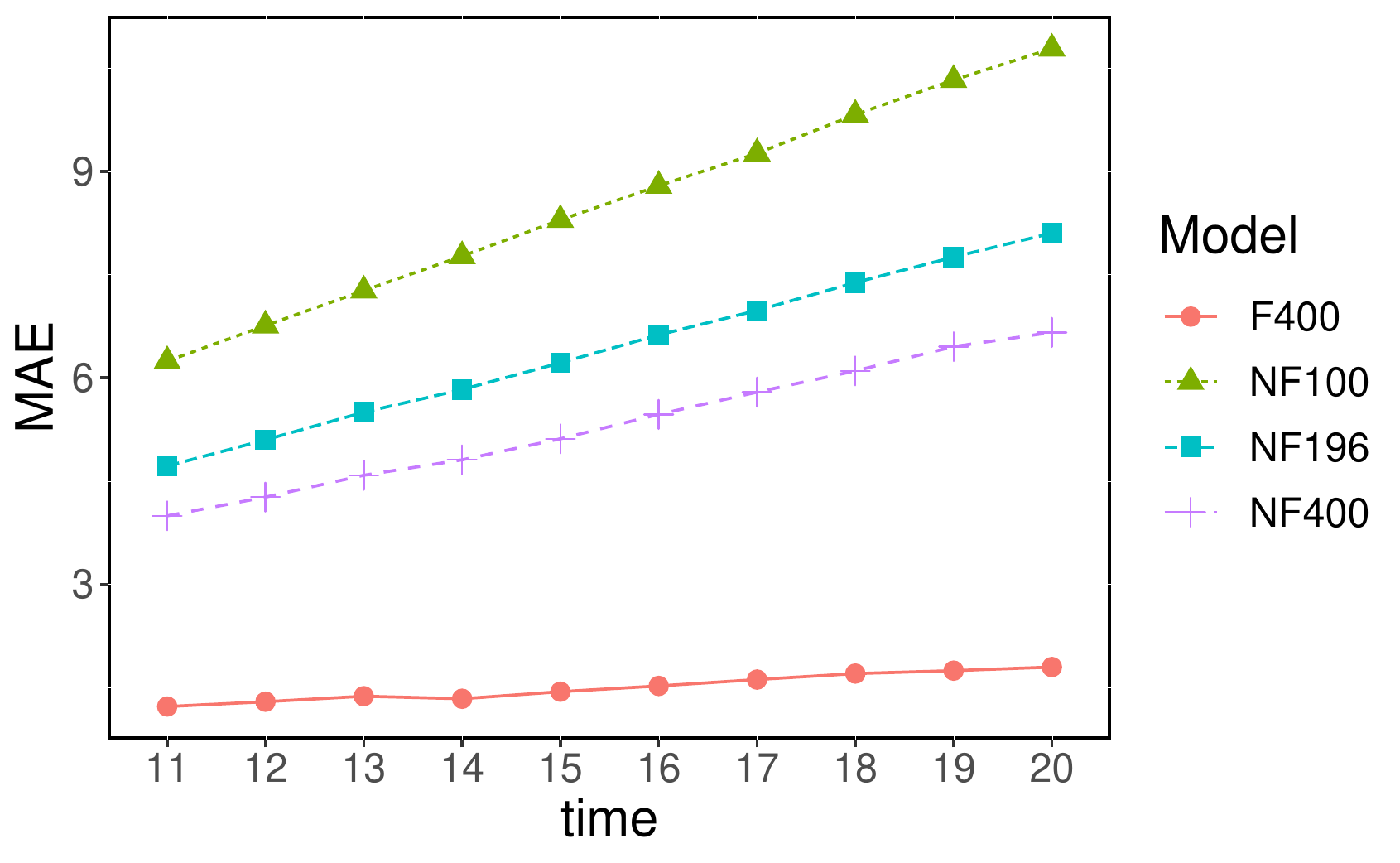}
    \caption{Comparison of MAE between four different models: the proposed approach effectively suppresses the Gibbs phenomenon and significantly improves the modeling accuracy.}\label{fig:SPcomparison}
\end{figure}

Moreover, because the Gibbs phenomenon occurs near the top boundary, Figure \ref{fig:SPcomparison} shows the MAE of the four models within the spatial area of $[0, 0.99]\times[0.95, 0.99]$ from time step 11 to time step 20. Note that, if there is no Gibbs phenomenon, the MAE over this spatial area is expected to be very close to 0 (as the strong signal is only found near the bottom boundary of the spatial domain). However, Figure \ref{fig:SPcomparison} clearly shows that the MAE associated with the existing approach is much higher than zero due to the Gibbs phenomenon (also see Figure \ref{fig:EstimatedSP}), while the MAE associated with the proposed approach is much smaller than that of the existing approach. Although keeping more high-frequency spectral coefficients does help the existing approach reduce the Gibbs phenomenon, \textit{the proposed approach is much more effective in terms of suppressing the Gibbs phenomenon, and as a result, significantly improves the modeling accuracy}. 

\emph{Discussion on the Hamming Window Method.} As mentioned in \citet{Pan1993}, a conventional way to suppress the Gibbs phenomenon is to use Hamming window to force the end points to become zero or near zero. By adopting the windows to the original signals, the discontinuity of the weighted signals is alleviated, and thus the Gibbs phenomenon can be suppressed to a certain extent. Here, we provide further discussions on the Hamming window method through the simulated data stream from Example I. 

\begin{figure}[h!]
    \centering
    \includegraphics[width=0.55\textwidth]{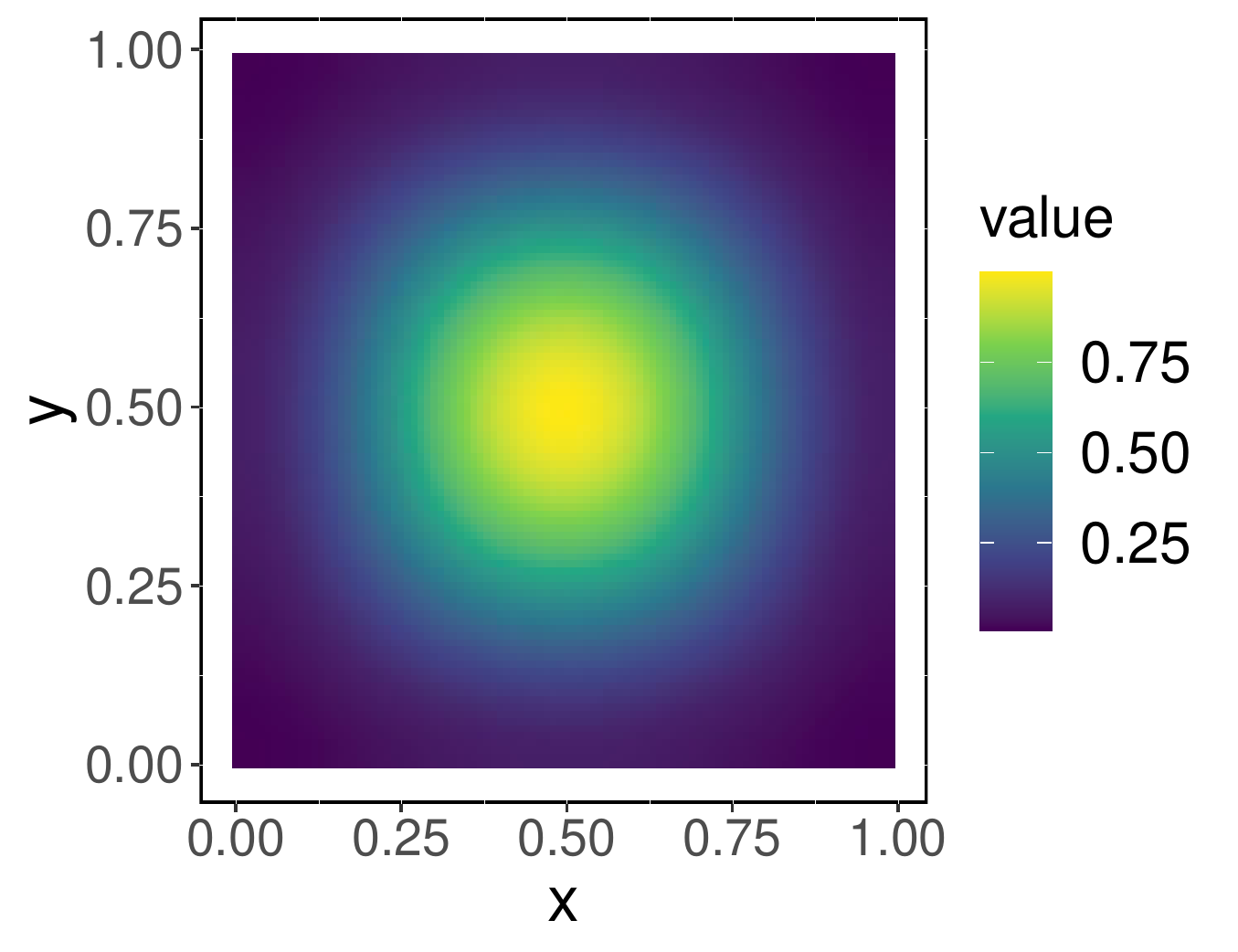}
    \caption{Hamming window for the simulated data stream}\label{fig:HammingWindow}
\end{figure}

For the Hamming window method, the 2-dimensional Hamming window can be defined as \citep{soon2003speech}:

\begin{equation}
    \bm{h} = \left(0.54-0.46\cos\left(\frac{2\pi \bm{i}}{99}\right)\right)\left(0.54-0.46\cos\left(\frac{2\pi \bm{j}^T}{99}\right)\right),
\end{equation}
where $\bm{i}=\bm{j}=(0,1,2,\cdots,99)^T$. Figure \ref{fig:HammingWindow} shows the 2-D Hamming window for the simulated data stream. We see that the boundaries are nearly zero. If the Hamming window is multiplied to $\bm{y}(t)$, i.e., $\text{vec}(\bm{h})\bm{y}(t)$, the Gibbs phenomenon can be suppressed due to the alleviated discontinuity in $\text{vec}(\bm{h})\bm{y}(t)$. 

However, adding the Hamming window causes the bias to the original data stream. Here, we first multiply the 2-D Hamming window to the simulated data stream and use the existing approach (\ref{eq:CaseDynamicalModel}) to model the transformed data stream. Then, we compare the modeling accuracy of 2-D hamming window method with our proposed flipping method. Again, we construct three models using the Hamming window method by letting $K=100, 196, 400$ in the existing approach (\ref{eq:CaseDynamicalModel}), and compare the MAEs of the three models and that of the proposed approach (\ref{eq:FSPDynamical}) with $K^* = 400$. These four models are denoted by ``HWNF100'', ``HWNF196'', ``HWNF400'', ``F400'', where ``HWNF'' means the Hamming window method with no process flipping, and ``F'' means the proposed approach with process flipping. 

Table \ref{tab:MAEComparison} presents the MAE of the four models over the entire spatial domain from time step 15 to time step 20. Note that, the modeling accuracy of the three models based on the Hamming window approach is not sensitive to the value of $K$. This is because the bias of these three models is mainly caused by introducing the Hamming window function to the original signals. Because the simulated data stream has large values on the bottom edge, multiplying the Hamming window inevitably leads to a high level of bias regardless of the value of $K$. We also see that all three models based on the Hamming window approach have higher MAEs than that of the proposed approach, suggesting that the Hamming window approach may not be the ideal choice for suppressing the Gibbs phenomenon in the statistical spatio-temporal model considered in this paper.

\begin{table}[h!]
\centering
\caption{Comparison of MAE between three models based on the Hamming window method and the model based on our the proposed approach}\label{tab:MAEComparison}
\begin{tabular}{c|cccccccccc} 
          \Xhline{1pt}
          Models & time 15& time 16& time 17& time 18& time 19& time 20\\
          \Xhline{1pt}
           HWNF100 & 3.694& 3.885& 4.091& 4.281& 4.466& 4.668\\
           HWNF196 & 3.686& 3.874& 4.082& 4.272& 4.458& 4.658\\
           HWNF400 & 3.679& 3.864& 4.071& 4.261& 4.446& 4.645\\
           F400 &  \textbf{1.681}& \textbf{1.745}& \textbf{1.814}& \textbf{1.861}& \textbf{1.923}& \textbf{1.995}\\
          \Xhline{1pt}
\end{tabular} 
\end{table}

\subsection{Example II: Modeling of a Precipitation Process}\label{sec:CaseStudy}
\begin{figure}[H]
    \centering
    \includegraphics[width=0.96\textwidth]{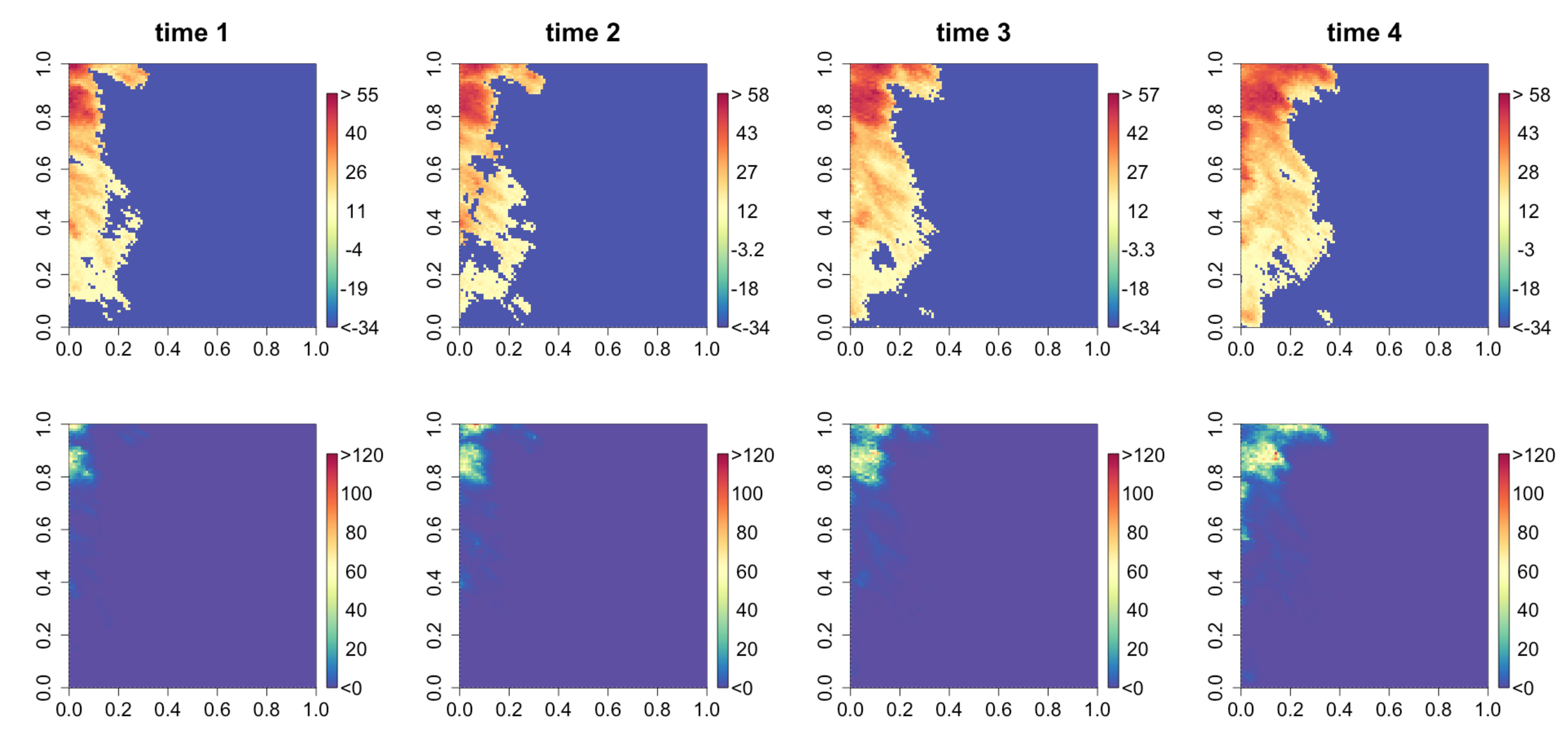}
    \vspace{-10pt}
    \caption{Weather radar images for a tropical storm; Top row: the original radar reflectivity image stream at 1 km above the mean sea level; Bottom row: the precipitation rate obtained from the radar reflectivity using the Marshall–Palmer relationship.}\label{fig:CaseRadar}
\end{figure}

In Example II, we apply the proposed method to a real dataset and demonstrate its capability of suppressing the Gibbs phenomenon. Figure \ref{fig:CaseRadar} shows a tropical precipitation process observed by a dual-polarization Meteorological Doppler Weather Radar (MDWR) system. The top row shows the standard Constant Altitude Plan Position Indicator (CAPPI) radar reflectivity data (in dBZ) at 1 km above the mean sea level, while the bottom row converts CAPPI data to precipitation rate (in mm/hr) using the empirical Marshall-Palmer relationship \citep{marshall1948distribution}:
\begin{equation}
    R = \left(\frac{10^{\frac{Z}{10}}}{200}\right)^{\frac{5}{8}}
\end{equation}
where $Z$ is the radar reflectivity in dBZ, and $R$ is the precipitation rate in mm$\cdot$hr$^{-1}$. Each image contains $100 \times 100$ pixels, and each pixel approximately covers a spatial area of $0.5\times0.5$ km$^2$. The time interval between two consecutive images is 5 mins.

\begin{figure}[h!]
    \centering
    \includegraphics[width=1\textwidth]{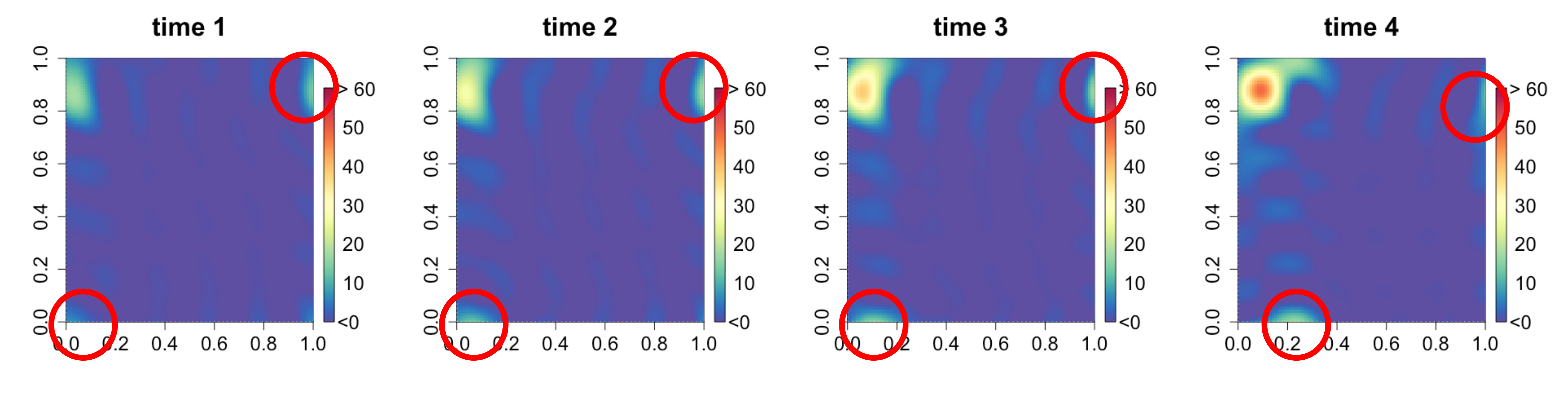}
    \caption{Both the Gibbs phenomenon and ripple artifact are clearly visible from the low-pass filtered precipitation process reconstructed by truncated Fourier series.}\label{fig:CaseLPRadar}
\end{figure}

In Figure \ref{fig:CaseRadar}, it is seen that the tropical thunderstorm system enters the spatial domain from the Northwest corner. Hence, the process is highly non-periodic and the images have boundary discontinuities. Such a non-periodic process gives rise to a strong Gibbs phenomenon. Figure \ref{fig:CaseLPRadar} shows the low-pass filtered precipitation process reconstructed by truncated the Fourier series. Both the Gibbs phenomenon and ripple artifact are clearly visible. The magnitude of the Gibbs phenomenon (in the South and East areas of the spatial domain) is close to the actual precipitation rate in the Northwest of the spatial domain, leading to erroneous precipitation prediction in the South and the East areas. Hence, for such a highly non-periodic precipitation process, suppressing the Gibbs phenomenon becomes critical in improving the accuracy of spatio-temporal modeling and prediction. 

\begin{figure}[h!]
    \centering
    \includegraphics[width=1\textwidth]{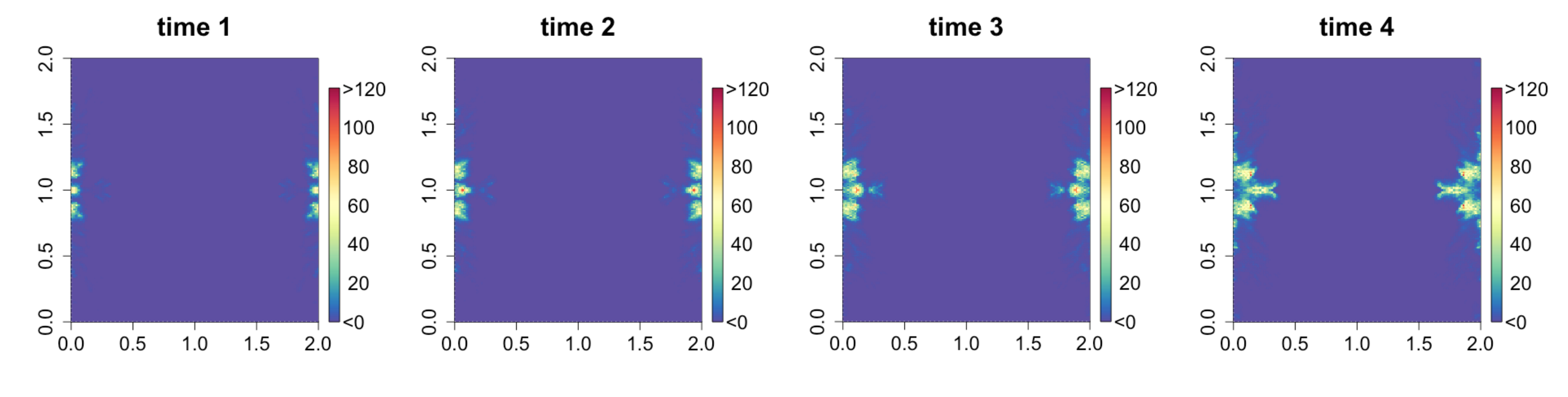}
    \caption{The flipped rainfall images are periodic with no boundary discontinuities}\label{fig:CaseFLPRadar}
\end{figure}

Following Section \ref{fig:flip}, we flip the radar images to make the spatial process periodic at each observation time. Here, $N_1=N_2=100$ and the flipped image has a total spatial points of $2N_1\times 2N_2=400$. Figure \ref{fig:CaseFLPRadar} shows the flipped radar images at times 1 to 4. 

\begin{figure}[ht]
    \centering
    \includegraphics[width=0.4\textwidth]{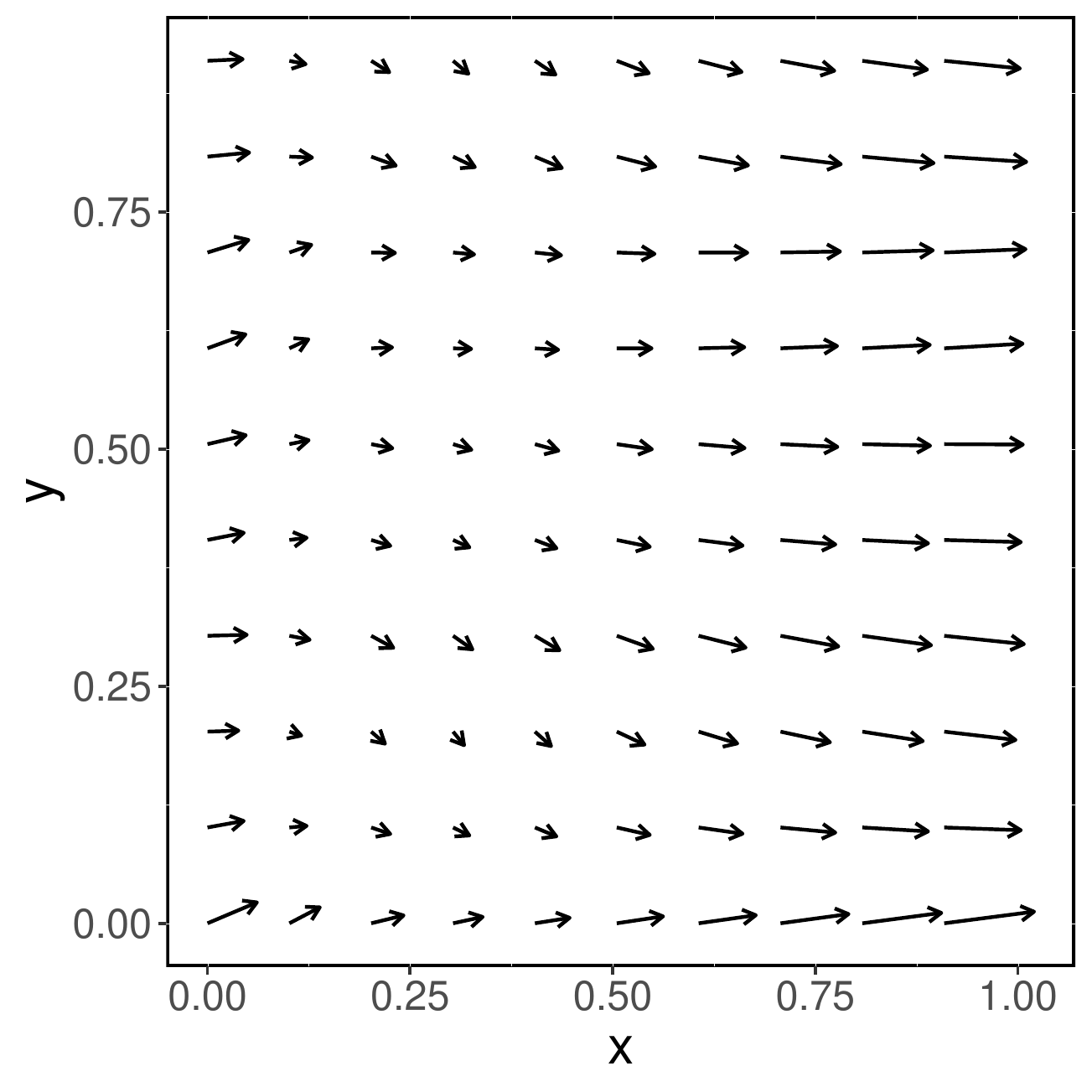}
    \caption{The velocity field obtained from the COTREC approach}\label{fig:CaseVelocity}
\end{figure}

The model (\ref{eq:CaseDynamicalModel}) can be constructed based on the flipped images. Because the transition matrix $\bm{G}^*(t)$ requires the information about the velocity field $\bm{v}(\bm{s}, t)$ and diffusivity fields $\bm{D}(\bm{s}, t)$, we adopt the COTREC (Tracking Radar Echoes by Correlation) approach to obtain the required velocity field from the radar reflectivity image stream \citep{li2004short}. COTREC is a widely used pattern-based method for tracking motion vectors of embedded radar echoes in the meteorological community. By partitioning the radar image into several overlapping small areas, the velocity of each area center can be obtained by maximizing the cross-correlation between that region and its surrounding areas from the consecutive images. Figure \ref{fig:CaseVelocity} shows the obtained velocity field from the COTREC approach, where the velocity speed ranges from 0 to 0.04 per time step. After the velocity field has been estimated, the diffusivity field is computed by $\bm{D}(\bm{s})=0.28(\delta_x\delta_y)\sqrt{\left(\partial v_x\big/\partial x-\partial v_y\big/\partial y\right)^2+\left(\partial v_x\big/\partial y+\partial v_y\big/\partial x\right)^2}$ where $\delta_x$ and $\delta_y$ are the computational resolution of the velocity field \citep{byun2006review, liu2018spatio}.

In addition, the likelihood function of $\sigma^2_{\bm{\alpha}}$ and $\sigma^2_{\bm{\beta}}$ can be obtained from the one-step-ahead predictive distribution of $\tilde{\bm{y}}^*(t)$ in (\ref{eq:CaseDynamicalModel}), and the maximum likelihood estimates of $\sigma^2_{\bm{\alpha}}$ and $\sigma^2_{\bm{\beta}}$ can thus be obtained. In this example, we use the \texttt{MARSS} package in R to obtain $\hat{\sigma}^2_{\bm{\alpha}}$ and $\hat{\sigma}^2_{\bm{\beta}}$. 

\begin{figure}[h!]
    \centering
    \includegraphics[width=1\textwidth]{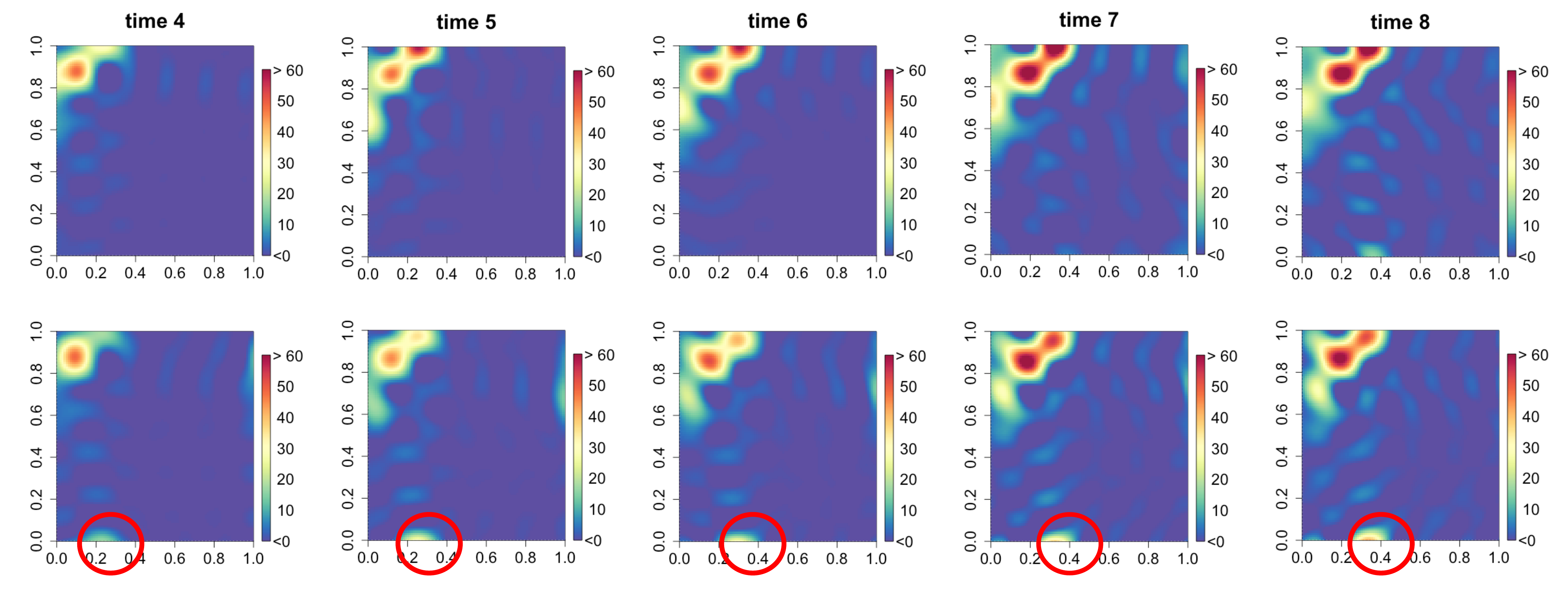}
    \caption{Comparison of the filtered and predicted images using the proposed approach (top) and the existing method (bottom): the proposed approach completely removes the Gibbs phenomena and significantly reduces the ripple artifact.}\label{fig:CaseFLTPRD}
\end{figure}

After the transition matrix $\bm{G}^*$, $\hat{\sigma}^2_{\bm{\alpha}}$ and $\hat{\sigma}^2_{\bm{\beta}}$ are obtained,  the state vector $\bm{\theta}(t)$ in (\ref{eq:CaseDynamicalModel}) is estimated by the Kalman Filter. In this example, we use 6 flipped radar images collected over a 25-minute interval to construct the model, and retain $K^*=200$ low-frequency spectral coefficients. Based on the filtered state variables  $\hat{\bm{\theta}}(t)$ and the Inverse Fourier Transform of the filtered states, the first row of Figure \ref{fig:CaseFLTPRD} shows the last three filtered and first two predicted images using the proposed approach. 
For comparison purposes, the second row of Figure \ref{fig:CaseFLTPRD} also presents the filtered and predicted images using the existing methods without image flipping. By comparing the results, we immediately see that \textit{the proposed approach completely removes the Gibbs phenomena and significantly reduces the ripple artifact}, potentially leading to more accurate precipitation forecast and producing higher quality precipitation maps. 

\section{Conclusions}
This paper provided an effective solution to suppress the Gibbs phenomenon in modeling PDE-based spatial-temporal processes due to boundary discontinuities. We proposed to flip the original spatial process at each time point, without altering the governing physics of the underlying process, to create a new periodic process that completely removes the discontinuity at boundaries. For the flipped process, we obtained the temporal dynamics of the spectral coefficients (associated with the flipped process), and integrated such dynamics into a dynamical model by extending a recently proposed physics-informed spatio-temporal model. Because the flipped spatio-temporal process is periodic and has a complete waveform, truncating the Fourier series on the flipped process no longer gives rise to the Gibbs phenomenon. Numerical investigations based on both simulated data and the radar image data of a tropical thunderstorm have been performed. The results obtained from the numerical investigations strongly demonstrate that the proposed approach successfully removes the Gibbs phenomenon, significantly reduces the ripple artifact due to truncation of the Fourier series, and potentially improves the modeling and prediction accuracy of statistical spatio-temporal models. 
Computer code is available at \url{https://github.com/gz-wei/Gibbs-Phenomena-Suppression-in-PDE-Based-Statistical-Spatio-Temporal-Models}. 

\section*{Appendix A}
$\bm{\Psi}^{(\text{A1+D1})}_{11}$ is a $|\Omega_1|\times|\Omega_1|$ matrix with the $(i,j)$-th,
element being given by $\Psi^{(\text{A1})}(\bm{k}_{1,j},\bm{k}_{1,i})+ \Psi^{(\text{D1})}(\bm{k}_{1,j},\bm{k}_{1,i})$, $\bm{\Psi}^{(\text{A1+D1})}_{12}$ is a $|\Omega_1|\times|\Omega_2|$ matrix with the $(i,j)$-th 
element being given by $\Psi^{(\text{A1})}(\bm{k}_{1,j},\bm{k}_{2,i})+ \Psi^{(\text{D1})}(\bm{k}_{1,j},\bm{k}_{2,i})$, $\bm{\Psi}^{(\text{A1+D1})}_{21}$ is a $|\Omega_2|\times|\Omega_1|$ matrix with the $(i,j)$-th,
element being given by $\Psi^{(\text{A1})}(\bm{k}_{1,j},\bm{k}_{2,i})+ \Psi^{(\text{D1})}(\bm{k}_{1,j},\bm{k}_{2,i})$, $\bm{\Psi}^{(\text{A1+D1})}_{22}$ is a $|\Omega_2|\times|\Omega_2|$ matrix with the $(i,j)$-th
element being given by $\Psi^{(\text{A1})}(\bm{k}_{2,j},\bm{k}_{2,i})+ \Psi^{(\text{D1})}(\bm{k}_{2,j},\bm{k}_{2,i})$, and 
\begin{align*}
    \Psi^{(\text{A1})}(\bm{k},\bm{k}') & =\int_{\mathbb{S}} \bm{v}^T(\bm{s},t)\tilde{\bm{k}}f_{\bm{k}}^{(s)}(\bm{s})f_{\bm{k}'}^{(c)}(\bm{s})d\bm{s}\\
    \Psi^{(\text{D1})}(\bm{k}, \bm{k}') & =\int_{\mathbb{S}}\left(-\tilde{\bm{k}}^T\bm{D}(\bm{s},t)\tilde{\bm{k}}f_{\bm{k}}^{(c)}(\bm{s}) -[\nabla\cdot\bm{D}(\bm{s}, t)]^T\tilde{\bm{k}}f_{\bm{k}}^{(s)}(\bm{s})\right)f^{(c)}_{\bm{k}'}(\bm{s})d\bm{s}
\end{align*}
with $\tilde{\bm{k}}=2\pi\bm{k}$.

$\bm{\Psi}^{(\text{A2+D2})}_{12}$ is a $|\Omega_1|\times|\Omega_2|$ matrix with the $(i,j)$-th element being given by $\Psi^{(\text{A2})}(\bm{k}_{1,j},\bm{k}_{2,i})+\Psi^{(\text{D2})}(\bm{k}_{1,j},\bm{k}_{2,i})$, $\bm{\Psi}^{(\text{A2+D2})}_{22}$ is a $|\Omega_2|\times|\Omega_2|$ matrix with the $(i,j)$-th element being given by $\Psi^{(\text{A2})}(\bm{k}_{2,j},\bm{k}_{2,i})+\Psi^{(\text{D2})}(\bm{k}_{2,j},\bm{k}_{2,i})$, and 
\begin{align*}
    \Psi^{(\text{A2})}(\bm{k},\bm{k})& =\int_{\mathbb{S}} -\bm{v}^T(\bm{s},t)\tilde{\bm{k}}f_{\bm{k}}^{(c)}(\bm{s})f_{\bm{k}'}^{(c)}(\bm{s})d\bm{s}\\
    \Psi^{(\text{D2})}(\bm{k}, \bm{k}')& =\int_{\mathbb{S}}\left(-\tilde{\bm{k}}^T\bm{D}(\bm{s}, t)\tilde{\bm{k}}f_{\bm{k}}^{(s)}(\bm{s})-[\nabla\cdot\bm{D}(\bm{s},t)]^T\tilde{\bm{k}}f_{\bm{k}}^{(c)}(\bm{s})\right)f^{(c)}_{\bm{k}'}(\bm{s})d\bm{s}.
\end{align*}

$\bm{\Psi}^{(\text{A3+D3})}_{21}$ is a $|\Omega_2|\times|\Omega_1|$ matrix with the $(i,j)$-th element being given by $\Psi^{(\text{A3})}(\bm{k}_{1,j},\bm{k}_{2,i})+\Psi^{(\text{D3})}(\bm{k}_{1,j},\bm{k}_{2,i})$, $\bm{\Psi}^{(\text{A3+D3})}_{22}$ is a $|\Omega_2|\times|\Omega_2|$ matrix with the $(i,j)$-th element being given by $\Psi^{(\text{A3})}(\bm{k}_{2,j},\bm{k}_{2,i})+\Psi^{(\text{D3})}(\bm{k}_{2,j},\bm{k}_{2,i})$, and 
\begin{align*}
   \Psi^{(\text{A3})}(\bm{k},\bm{k}')& = \int_{\mathbb{S}} \bm{v}^T(\bm{s},t)\tilde{\bm{k}}f_{\bm{k}}^{(s)}(\bm{s})f_{\bm{k}'}^{(s)}(\bm{s})d\bm{s}\\
    \Psi^{(\text{D3})}(\bm{k}, \bm{k}') & = \int_{\mathbb{S}}\left(-\widetilde{\bm{k}}^T\bm{D}(\bm{s}, t)\widetilde{\bm{k}}f_{\bm{k}}^{(c)}(\bm{s})-[\nabla\cdot\bm{D}(\bm{s},t)]^T\widetilde{\bm{k}}f_{\bm{k}}^{(s)}(\bm{s})\right)f^{(s)}_{\bm{k}'}(\bm{s})d\bm{s}.
\end{align*}

$\bm{\Psi}^{(\text{A4+D4})}_{22}$ is a $|\Omega_2|\times|\Omega_2|$ matrix with the $(i,j)$-th element being given by $\Psi^{(\text{A4})}(\bm{k}_{2,j},\bm{k}_{2,i})+\Psi^{(\text{D4})}(\bm{k}_{2,j},\bm{k}_{2,i})$, and 
\begin{align*}
    \Psi^{(\text{A4})}(\bm{k},\bm{k}')& =-\int_{\mathbb{S}} \bm{v}^T(\bm{s},t)\tilde{\bm{k}}f_{\bm{k}}^{(c)}(\bm{s})f_{\bm{k}'}^{(s)}(\bm{s})d\bm{s}\\
    \Psi^{(\text{D4})}(\bm{k}, \bm{k}')& =\int_{\mathbb{S}}\left(-\tilde{\bm{k}}^T\bm{D}(\bm{s}, t)\tilde{\bm{k}}f_{\bm{k}}^{(s)}(\bm{s})-[\nabla\cdot\bm{D}(\bm{s},t)]^T\tilde{\bm{k}}f_{\bm{k}}^{(c)}(\bm{s})\right)f^{(s)}_{\bm{k}'}(\bm{s})d\bm{s}.
\end{align*}

$\bm{C}_1$ is a diagonal matrix $\text{diag}( \{c_{\bm{k}_{1,i}} \}_{i=1}^{|\mathcal{K}_1|})$, and $\bm{C}_2$ is a diagonal matrix $\text{diag}( \{c_{\bm{k}_{2,i}} \}_{i=1}^{|\mathcal{K}_2|})$ where
\begin{equation}
	c_{\bm{k}}=\int_{\mathbb{S}}\cos^2(\bm{k}^T\bm{s})d\bm{s}.
\end{equation}


\begin{thebibliography}{}
\bibitem[{Guinness and Stein, 2013}]{guinness2013interpolation}
Guinness, J. and Stein, M. L. (2013). Interpolation of non-stationary high frequency spatial- temporal temperature. \emph{The Annals of Applied Statistics}, pages 1684–1708.

\bibitem[{Kuusela and Stein, 2018}]{kuusela2018locally}
Kuusela, M. and Stein, M. L. (2018). Locally stationary spatio-temporal interpolation of argo profiling float data. \emph{Proceedings of the Royal Society A}, 474(2220):20180400.

\bibitem[{Vandeskog et al., 2013}]{vandeskog2022quantile}
Vandeskog, S. M., Thorarinsdottir, T. L., Steinsland, I., and Lindgren, F. (2022). Quantile based modeling of diurnal temperature range with the five-parameter lambda distribution. \emph{Environmetrics}, 33(4):e2719. \href{https://doi.org/10.1002/env.2719}{https://doi.org/10.1002/env.2719}  

\bibitem[{Huang and Hsu, 2004}]{huang2004modeling}
Huang, H.-C. and Hsu, N.-J. (2004). Modeling transport effects on ground-level ozone using a non-stationary space-time model. \emph{Environmetrics}, 15(3):251–268.

\bibitem[{Schliep et al., 2020}]{schliep2020data}
Schliep, E. M., Collins, S. M., Rojas-Salazar, S., Lottig, N. R., and Stanley, E. H. (2020). Data fusion model for speciated nitrogen to identify environmental drivers and improve estimation of nitrogen in lakes. \emph{The Annals of Applied Statistics}, 14(4):1651–1675.

\bibitem[{Liu, et al., 2016}]{liu2016statistical}
Liu, X., Yeo, K., Hwang, Y., Singh, J., and Kalagnanam, J. (2016). A statistical modeling approach for air quality data based on physical dispersion processes and its application to ozone modeling. \emph{The Annals of Applied Statistics}, 10(2):756–785.

\bibitem[{Schliep et al., 2020}]
Schliep, E. M., Collins, S. M., Rojas-Salazar, S., Lottig, N. R., and Stanley, E. H. (2020). Data fusion model for speciated nitrogen to identify environmental drivers and improve estimation of nitrogen in lakes. \emph{The Annals of Applied Statistics}, 14(4):1651–1675.

\bibitem[{Fioravanti et al., 2022}]{fioravanti2022spatiotemporal}
Fioravanti, G., Cameletti, M., Martino, S., Cattani, G., and Pisoni, E. (2022). A spatiotemporal analysis of no2 concentrations during the italian 2020 covid-19 lockdown. \emph{Environmetrics}, page e2723. \href{https://doi.org/10.1002/env.2723}{https://doi.org/10.1002/env.2723} 

\bibitem[{Cooley et al., 2007}]{cooley2007bayesian}
Cooley, D., Nychka, D., and Naveau, P. (2007). Bayesian spatial modeling of extreme precipitation return levels. \emph{Journal of the American Statistical Association}, 102(479):824– 840.

\bibitem[{Heaton et al., 2011}]{heaton2011spatio}
Heaton, M. J., Katzfuss, M., Ramachandar, S., Pedings, K., Gilleland, E., Mannshardt- Shamseldin, E., and Smith, R. L. (2011). Spatio-temporal models for large-scale indicators of extreme weather. \emph{Environmetrics}, 22(3):294–303.

\bibitem[{Kleiber et al., 2012}]{kleiber2012daily}
Kleiber, W., Katz, R. W., and Rajagopalan, B. (2012). Daily spatio-temporal precipitation simulation using latent and transformed gaussian processes. Water Resources Research, 48(1).

\bibitem[{Liu et al., 2018}]{liu2018spatio}
Liu, X., Gopal, V., and Kalagnanam, J. (2018). A spatio-temporal modeling framework for weather radar image data in tropical southeast asia. \emph{The Annals of Applied Statistics}, 12(1):378–407.

\bibitem[{Bopp et al., 2021}]{bopp2021hierarchical}
Bopp, G. P., Shaby, B. A., and Huser, R. (2021). A hierarchical max-infinitely divisible spatial model for extreme precipitation. \emph{Journal of the American Statistical Association}, 116(533):93–106.

\bibitem[{Wei et al., 2022}]{wei2022}
Wei, G., Krishnan, V., Xie, Y., Sengupta, M., Zhang, Y., Liao, H., and Liu, X. (2022). Physics-informed statistical modeling for wildfire aerosols process using multi- source geostationary satellite remote-sensing data streams. \emph{arXiv preprint arXiv:2206.11766}. \href{https://doi.org/10.48550/arXiv.2206.11766}{https://doi.org/10.48550/arXiv.2206.11766}

\bibitem[{Kang et al., 2011}]{kang2011meta}
Kang, J., Johnson, T. D., Nichols, T. E., and Wager, T. D. (2011). Meta analysis of functional neuroimaging data via bayesian spatial point processes. \emph{Journal of the American Statistical Association}, 106(493):124–134.

\bibitem[{Katzfuss and Cressie, 2011}]{katzfuss2011spatio}
Katzfuss, M. and Cressie, N. (2011). Spatio-temporal smoothing and em estimation for massive remote-sensing data sets. \emph{Journal of Time Series Analysis}, 32(4):430–446.

\bibitem[{Hefley et al., 2017}]{hefley2017dynamic}
Hefley, T. J., Hooten, M. B., Hanks, E. M., Russell, R. E., and Walsh, D. P. (2017). Dynamic spatio-temporal models for spatial data. \emph{Spatial Statistics}, 20:206–220.

\bibitem[{Reich et al., 2018}]{reich2018fully}
Reich, B. J., Guinness, J., Vandekar, S. N., Shinohara, R. T., and Staicu, A.-M. (2018). Fully bayesian spectral methods for imaging data. \emph{Biometrics}, 74(2):645–652.

\bibitem[{Castruccio et al., 2018}]{castruccio2018scalable}
Castruccio, S., Ombao, H., and Genton, M. G. (2018). A scalable multi-resolution spatio- temporal model for brain activation and connectivity in fmri data. \emph{Biometrics}, 74(3):823– 833.

\bibitem[{Zammit-Mangion et al., 2021}]{zammit2021deep}
Zammit-Mangion, A., Ng, T. L. J., Vu, Q., and Filippone, M. (2021). Deep compositional spatial models. \emph{Journal of the American Statistical Association}, pages 1–22.

\bibitem[{Zimmerman, 2006}]{zimmerman2006optimal}
Zimmerman, D. L. (2006). Optimal network design for spatial prediction, covariance param- eter estimation, and empirical prediction. \emph{Environmetrics}, 17(6):635–652.

\bibitem[{Zimmerman and Buckland, 2019}]{zimmerman2019environmental}
Zimmerman, D. and Buckland, S. (2019). \emph{``Environmental sampling design'', in handbook of environmental and ecological statistics}. CRC Press, pp. 181-210.

\bibitem[{Sharrock and Kantas, 2022}]{sharrock2022joint}
Sharrock, L. and Kantas, N. (2022). Joint online parameter estimation and optimal sensor placement for the partially observed stochastic advection-diffusion equation. \emph{SIAM/ASA Journal on Uncertainty Quantification}, 10(1):55–95.

\bibitem[{Cressie and Wikle, 2011}]{Cressie2011}
Cressie, N. and Wikle, C. (2011). \emph{Statistics for spatio-temporal data}. John Wiley \& Sons, Hoboken, New Jersey.

\bibitem[{Sigrist et al., 2015}]{Sigrist2015}
Sigrist, F., Kunsch, H. R., and Stahel, W. A. (2015). Stochastic partial differential equation based modelling of large space-time data sets. \emph{Journal of the Royal Statistical Society: Series B}, 77(1):3–33.

\bibitem[{Kutz et al., 2016}]{Kutz2016}
Kutz, N., Brunton, S. L., W., B. B., and Proctor, J. L. (2016). \emph{Dynamic mode decomposition: data-driven modeling of complex systems}. SIAM-Society for Industrial and Applied Mathematics.

\bibitem[{Liu et al., 2021}]{Liu2021}
Liu, X., Yeo, K., and Lu, S. (2021). Statistical modeling for spatio-temporal data from stochastic convection-diffusion processes. \emph{Journal of the American Statistical Association}. \href{https://doi.org/10.1080/01621459.2020.1863223}{https://doi.org/10.1080/01621459.2020.1863223} 

\bibitem[{Pan, 1993}]{Pan1993}
Pan, C. (1993). \emph{Gibbs phenomenon suppression and optimal windowing for attenuation and Q measurements}. SLAC National Accelerator Laboratory: SLAC-PUB-6222.

\bibitem[{Jerri, 1998}]{Jerri1998}
Jerri, A. J. (1998). \emph{The Gibbs phenomenon in fourier analysis}, Splines and Wavelet Approximations. Springer.

\bibitem[{Prabhu, 2014}]{Prabhu2014}
Prabhu, K. M. M. (2014). \emph{Window functions and their applications in signal processing}. CRC Press.

\bibitem[{Pan, 2001}]{Pan2001}
Pan, C. (2001). Gibbs phenomenon removal and digital filtering directly through the fast fourier transform. \emph{IEEE Transactions on Signal Processing}, 49:444–447.

\bibitem[{Wikle et al., 1998}]{wikle1998hierarchical}
Wikle, C. K., Berliner, L. M., and Cressie, N. (1998). Hierarchical bayesian space-time models. \emph{Environmental and Ecological Statistics}, 5(2):117–154.

\bibitem[{Evensen and Van Leeuwen, 2000}]{evensen2000ensemble}
Evensen, G. and Van Leeuwen, P. J. (2000). An ensemble kalman smoother for nonlinear dynamics. \emph{Monthly Weather Review}, 128(6):1852–1867.

\bibitem[{Stroud and Bengtsson, 2007}]{stroud2007sequential}
Stroud, J. R. and Bengtsson, T. (2007). Sequential state and variance estimation within the ensemble kalman filter. \emph{Monthly Weather Review}, 135(9):3194–3208.

\bibitem[{Stroud et al., 2010}]{stroud2010ensemble}
Stroud, J. R., Stein, M. L., Lesht, B. M., Schwab, D. J., and Beletsky, D. (2010). An ensemble kalman filter and smoother for satellite data assimilation. \emph{Journal of the American Statistical Association}, 105(491):978–990.

\bibitem[{Soon and Koh, 2003}]{soon2003speech}
Soon, Y. and Koh, S. N. (2003). Speech enhancement using 2-d fourier transform. \emph{IEEE Transactions on Speech and Audio Processing}, 11(6):717–724.

\bibitem[{Marshall, 1948}]{marshall1948distribution}
Marshall, J. S. (1948). The distribution of raindrops with size. \emph{Journal of Meteorology}., 5:165–166.

\bibitem[{Li and Lai, 2004}]{li2004short}
Li, P. and Lai, E. S. (2004). Short-range quantitative precipitation forecasting in hong kong. \emph{Journal of Hydrology}, 288(1-2):189–209.

\bibitem[{Byun and Schere, 2006}]{byun2006review}
Byun, D. and Schere, K. L. (2006). Review of the governing equations, computational algorithms, and other components of the models-3 community multiscale air quality (cmaq) modeling system. \emph{Applied Mechanics Review}, (2006): 51-77

\end{thebibliography}
\end{document}